\def\abstracttext{Consider a diagram of quasi-categories that admit and functors that preserve limits or colimits of a fixed shape. We show that any weighted limit whose weight is a projective cofibrant simplicial functor is again a quasi-category admitting these (co)limits and that they are preserved by the functors in the limit cone. In particular, the Bousfield-Kan homotopy limit of a diagram of quasi-categories admits any (co)limits existing in and preserved by the functors in that diagram. In previous work, we demonstrated that the quasi-category of algebras for a homotopy coherent monad could be described as a weighted limit of this type, so these results specialise to (co)completeness results for quasi-categories of algebras. The second half of this paper establishes a further result in the quasi-categorical setting: proving, in analogy with the classical categorical case, that the monadic forgetful functor of the quasi-category of algebras for a homotopy coherent monad creates all limits that exist in the base quasi-category, regardless of whether its functor part preserves those limits. This proof relies upon a more delicate and explicit analysis of the particular weight used to define quasi-categories of algebras. }

\expandafter\ifx\csname c@nestcount\endcsname\relax\newcounter{nestcount}\fi
\global\addtocounter{nestcount}{1}
\ifnum\value{nestcount}>1\else

\documentclass[12pt,reqno]{amsart}

\usepackage{a4wide}
\usepackage{amsmath}
\usepackage{amssymb}
\usepackage{amsxtra}
\usepackage{amsthm}
\usepackage{mathtools}
\usepackage{extpfeil}
\usepackage{enumitem}
\usepackage[mathscr]{eucal}
\usepackage{graphicx}
\usepackage{asymptote}
\usepackage{ifnextok}
\usepackage[T1]{fontenc}

\usepackage[all]{xy}
\SelectTips{cm}{}
\SilentMatrices
\ifx\pdfoutput\undefined
  \xyoption{dvips}
\fi

\usepackage[silent]{fixme}
\usepackage[unicode,final]{hyperref}
\usepackage{slashed}

\usepackage{mathrsfs}

\newcommand{\ifundef}[1]{\expandafter\ifx\csname#1\endcsname\relax}

\DeclareMathAlphabet{\mathbbe}{U}{bbold}{m}{n}

\makeatletter

\def\re@DeclareMathSymbol#1#2#3#4{%
    \let#1=\undefined
    \DeclareMathSymbol{#1}{#2}{#3}{#4}}

\ifundef{Top}
  \DeclareSymbolFont{tcSyC}{U}{txsyc}{m}{n}
  \SetSymbolFont{tcSyC}{bold}{U}{txsyc}{bx}{n}
  \DeclareFontSubstitution{U}{txsyc}{m}{n}

  \re@DeclareMathSymbol{\Top}{\mathord}{tcSyC}{120}
  \re@DeclareMathSymbol{\Bot}{\mathord}{tcSyC}{121}
\fi

\ifundef{righthalfcup}
  \DeclareFontFamily{U}{MnSymbolC}{}
  \DeclareSymbolFont{mnSyC}{U}{MnSymbolC}{m}{n}
  \SetSymbolFont{mnSyC}{bold}{U}{MnSymbolC}{b}{n}
  \DeclareFontShape{U}{MnSymbolC}{m}{n}{
      <-6>  MnSymbolC5
     <6-7>  MnSymbolC6
     <7-8>  MnSymbolC7
     <8-9>  MnSymbolC8
     <9-10> MnSymbolC9
    <10-12> MnSymbolC10
    <12->   MnSymbolC12}{}
  \DeclareFontShape{U}{MnSymbolC}{b}{n}{
      <-6>  MnSymbolC-Bold5
     <6-7>  MnSymbolC-Bold6
     <7-8>  MnSymbolC-Bold7
     <8-9>  MnSymbolC-Bold8
     <9-10> MnSymbolC-Bold9
    <10-12> MnSymbolC-Bold10
    <12->   MnSymbolC-Bold12}{}
  
  \re@DeclareMathSymbol{\righthalfcup}{\mathord}{mnSyC}{184}
  \re@DeclareMathSymbol{\lefthalfcap}{\mathord}{mnSyC}{185}
\fi

\DeclareFontFamily{U}{MnSymbolA}{}
\DeclareSymbolFont{mnSyA}{U}{MnSymbolA}{m}{n}
\SetSymbolFont{mnSyA}{bold}{U}{MnSymbolA}{b}{n}
\DeclareFontShape{U}{MnSymbolA}{m}{n}{
    <-6>  MnSymbolA5
   <6-7>  MnSymbolA6
   <7-8>  MnSymbolA7
   <8-9>  MnSymbolA8
   <9-10> MnSymbolA9
  <10-12> MnSymbolA10
  <12->   MnSymbolA12}{}
\DeclareFontShape{U}{MnSymbolA}{b}{n}{
    <-6>  MnSymbolA-Bold5
   <6-7>  MnSymbolA-Bold6
   <7-8>  MnSymbolA-Bold7
   <8-9>  MnSymbolA-Bold8
   <9-10> MnSymbolA-Bold9
  <10-12> MnSymbolA-Bold10
  <12->   MnSymbolA-Bold12}{}

\re@DeclareMathSymbol{\twoheadedswarrow}{\mathord}{mnSyA}{30}

\makeatother

\newcommand{\mlaux}[3]{\setbox0=\hbox{$\mathsurround=0pt #2{#3}$}%
  \dimen0=\dp0\advance\dimen0 by \ht0\lower#1\dimen0\box0}

\newcommand{\makellapm}[2]{\hbox to 0pt{\hss$\mathsurround=0pt #1{#2}$}}

\newcommand{\makerlapm}[2]{\hbox to 0pt{$\mathsurround=0pt #1{#2}$\hss}}

\newcommand{\makelapm}[2]{\hbox to 0pt{\hss$\mathsurround=0pt #1{#2}$\hss}}

\newcommand{\makeushort}[3]{%
	\setbox0=\hbox{$\mathsurround=0pt #2{#3}$}%
	\hbox to 1\wd0{\hss\underbar{\hbox to #1\wd0{\hss\box0\hss}}\hss}}

\def\makebigger#1#2#3{\scalebox{#1}{$\mathsurround=0pt #2{#3}$}}
\def\bigger#1#2{{\relax\mathpalette{\makebigger{#1}}{#2}}}

\def\scaleuphalf{1.0954}
\def\scaleupone{1.2}

\newcommand{\op}{^{\mathord{\text{\rm op}}}}

\renewcommand{\th}{^{\text{th}}}

\newcommand{\defeq}{\mathrel{:=}}

\makeatletter

\def\newmop{\@ifstar{\@newmop m}{\@newmop o}}
\def\@newmop#1{\@ifnextchar[{\@@newmop #1}{\@@@newmop #1}}
\def\@@newmop#1[#2]{\@declmathop #1#2}
\def\@@@newmop#1#2{\expandafter\@declmathop\expandafter #1\csname #2\endcsname{#2}}

\makeatother

\newdir{ |}{{}*!/-5pt/@{|}}

\newmop{im}
\newmop{coim}
\newmop{dom}
\newmop{cod}
\newmop{id}

\newmop{obj}
\newmop{arr}
\newmop{sq}
\newmop{norm}

\newmop{el}

\newmop{ev}

\newmop{Ext}
\newmop{icon}

\newmop{cell}
\newmop{cof}
\newmop{fib}

\newmop{diag}

\newmop*{colim}
\newmop*{holim}

\newmop[\lan]{lan}
\newmop[\ran]{ran}

\newcommand{\comma}{\mathbin{\downarrow}}

\makeatletter
\newcommand{\rotatemath}[2]{\rotatebox[origin=c]{180}{$\m@th #1{#2}$}}
\makeatother

\newmop[\pr]{pr}
\newmop[\dm]{d} 

\newmop[\cpr]{in}

\newcommand{\pocorner}{\hbox to 8pt{{\vrule height8pt depth0pt width0.5pt}%
    \vbox to 8pt{{\hrule height0.5pt width7.5pt depth0pt}\vfill}}}
\newcommand{\poexcursion}{\save[]-<15pt,-15pt>*{\pocorner}\restore}
\newcommand{\pbcorner}{\vbox to 0pt{\kern 4pt\hbox to 0pt{\kern 4pt%
      \vbox{{\hrule height0.5pt width7.5pt depth0pt}}%
      {\vrule height8pt depth0pt width0.5pt}\hss}\vss}}
\newcommand{\pbexcursion}{\save[]+<5pt,-5pt>*{\pbcorner}\restore}

\newcommand{\wlim}[2]{\{#1,#2\}}

\newmop{cls}

\newcommand{\category}[1]{\underline{\smash[b]{\text{\rm{#1}}}}}

\newcommand{\tcat}{\lcat}

\newcommand{\lcat}{\mathcal}
\newcommand{\scat}{\mathbf}
\newcommand{\stcat}{\scat}

\newcommand{\catone}{{\bigger{1.16}{\mathbbe{1}}}}

\makeatletter

\def\Del@Sym{{\bigger\scaleuphalf{\mathbbe{\Delta}}}}

\def\del@fn{\futurelet\del@next}
\def\del@dn{\def\del@next}

\def\parsedel@{%
  \ifx +\del@next \del@dn+{\Del@Sym_{\mathord{+}}}%
  \else \del@dn {\del@fn\parsedel@@}%
  \fi\del@next}

\def\parsedel@@{%
  \ifx\space@\del@next \expandafter\del@dn\space{\del@fn\parsedel@@}%
  \else\ifx [\del@next \del@dn[{\del@fn\parsedel@@@}%
  \else\ifx _\del@next \del@dn{\Delta}%
  \else\ifx ^\del@next \del@dn{\Delta}%
  \else \del@dn{\Del@Sym}%
  \fi\fi\fi\fi\del@next}

\def\parsedel@@@{%
  \ifx\space@\del@next \expandafter\del@dn\space{\del@fn\parsedel@@@}%
  \else\ifx t\del@next \del@dn t{\Del@Sym_\infty\del@fn\parsedel@@@@}%
  \else\ifx b\del@next \del@dn b{\Del@Sym_{-\infty}\del@fn\parsedel@@@@}%
  \else \del@dn{\errmessage{unexpected modifier}}%
  \fi\fi\fi\del@next}

\def\parsedel@@@@{%
  \ifx\space@\del@next \expandafter\del@dn\space{\del@fn\parsedel@@@@}%
  \else\ifx ]\del@next \del@dn]{}%
  \else \del@dn{\errmessage{expecting close of option block}}%
  \fi\fi\del@next}

\def\Del{\del@fn\parsedel@}

\makeatother

\newcommand{\sSet}{\category{sSet}}
\newcommand{\qCat}{\category{qCat}}
\newcommand{\Adj}{\category{Adj}}
\newcommand{\Mnd}{\category{Mnd}}

\newcommand{\pbshape}{{\mathord{\bigger\scaleupone\righthalfcup}}}
\newcommand{\poshape}{{\mathord{\bigger\scaleupone\lefthalfcap}}}

\newcommand{\fbv}[1]{\{{#1}\}}

\newmop{dec}
\newmop{fatdec}
\newmop{slc}
\newmop{fatslc}

\newcommand{\slice}{/}

\newcommand{\slicer}[2]{{#1\mkern-1mu}_{{}\slice{#2}}}

\newmop{ir} % the interval representation functor
\newmop{incl}

\newcommand{\ho}{h}

\newcommand{\boundary}{\partial}

\def\reedyfilt#1_#2{#1_{\leq #2}}
\newmop{sk}
\newmop{cosk}
\newmop{res}

\newmop{map}
\newmop{Ho}

\newmop[\pth]{path}
\newmop{cyl}

\newmop[\const]{c}
  
\newmop{coll}
\newmop{wgt}

\newdir{ >}{{}*!/-7pt/@{>}}
\newdir{u(}{{}*!/-4pt/@^{(}}
\newdir{d(}{{}*!/-4pt/@_{(}}
\newdir{|>}{%
  !/4.5pt/@{|}*:(1,-.2)@^{>}*:(1,+.2)@_{>}*+@{}}

\makeatletter
\def\makeslashed#1#2#3#4#5{#1{\mathpalette{\sla@{#2}{#3}{#4}}{#5}}}

\def\@mathlower#1#2#3{\setbox0=\hbox{$\m@th#2#3$}\lower#1\ht0\box0}
\def\mathlower#1#2{\mathpalette{\@mathlower{#1}}{#2}}
\makeatother

\newcommand{\inc}{\hookrightarrow}

\newcommand{\tfib}{\twoheadrightarrow}

\makeatletter

\def\tens@fn{\futurelet\tens@next}
\def\tens@dn{\def\tens@nextcont}
\newtoks\tens@toks
\def\addtotens@toks#1{\tens@toks=\expandafter{\the\tens@toks#1}}

\def\parsetens@@{%
    \ifx\space@\tens@next \expandafter\tens@dn\space{\tens@fn\parsetens@@}%
    \else\ifx ^\tens@next \tens@dn ^##1{\parsetens@procsep^\addtotens@toks{##1}%
      \tens@fn\parsetens@@}%
    \else\ifx _\tens@next \tens@dn _##1{\parsetens@procsep_\addtotens@toks{##1}%
      \tens@fn\parsetens@@}%
    \else\tens@dn{\ifx *\tens@last \else\addtotens@toks\egroup\fi\the\tens@toks}%
    \fi\fi\fi\tens@nextcont}

\def\parsetens@procsep#1{%
  \ifx *\tens@last \addtotens@toks{#1}\addtotens@toks\bgroup%
  \else\ifx \tens@last\tens@next \addtotens@toks,%
  \else \addtotens@toks\egroup\addtotens@toks\bgroup%
    \addtotens@toks\egroup\addtotens@toks{#1}\addtotens@toks\bgroup% 
  \fi\fi\let\tens@last\tens@next}

\newcommand{\tn}[1]{\let\tens@last=*\tens@toks={#1}\tens@fn\parsetens@@}

\makeatother

\def\adjdisplay#1-|#2:#3->#4.{{%
    \xymatrix@R=0em@!C=2.5em{%
      *+[l]{#3} \ar@/_0.55pc/[rr]_-{#2} & {\bot} &
      *+[r]{#4}\ar@/_0.55pc/[ll]_-{#1}}}}

\def\adjdisplaytwo#1-|#2:#3->#4.{{% 
\xymatrix@=1.2em{
      {#3}\ar@/_1.5ex/[rr]_-{#2}^-{}="one"
      & & {#4}
      \ar@/_1.5ex/[ll]_-{#1}^-{}="two" 
      \ar@{}"one";"two"|{\bot}
    }}}

\def\tripleadjdisplay#1-|#2-|#3:#4->#5.{{%
\xymatrix@=2.4em{ 
{#4}\ar[r]|{#2} &
{#5} \ar@/_3ex/[l]_{#1}^{\bot} \ar@/^3ex/[l]_{\bot}^{#3}}
}}

\def\adjinline#1-|#2:#3->#4.{{#1}\dashv{#2}:#3\to #4}

\newcommand{\pent}[1]{
  \xybox{
    \POS (0,-15)*+{\a}="0", 
         (-14,-5)*+{\b}="1", 
         (-9,12)*+{\c}="2", 
         (9,12)*+{\d}="3", 
         (14,-5)*+{\e}="4"
    \POS"0" \ar "1"^{\labelstyle \ab}|{}="01"
    \POS"1" \ar "2"^{\labelstyle \bc}|{}="12"
    \POS"2" \ar "3"^{\labelstyle \cd}|{}="23"
    \POS"3" \ar "4"^{\labelstyle \de}|{}="34"
    \POS"0" \ar "4"_{\labelstyle \ae}|{}="04"
    \ifcase #1
    \POS"0" \ar "2"|{\labelstyle \ac}="02"
    \POS"0" \ar "3"|{\labelstyle \ad}="03"
    \POS"02";"1"**{}, ?(0.3) \ar@{=>} ?(0.7)^{\labelstyle \abc}
    \POS"03";"2"**{}, ?(0.25) \ar@{=>} ?(0.5)_{\labelstyle \acd}
    \POS"04";"3"**{}, ?(0.2) \ar@{=>} ?(0.4)_{\labelstyle \ade}
    \or
    \POS"1" \ar "3"|{\labelstyle \bd}="13"
    \POS"1" \ar "4"|{\labelstyle \be}="14"
    \POS"13";"2"**{}, ?(0.3) \ar@{=>} ?(0.7)_{\labelstyle \bcd}
    \POS"14";"3"**{}, ?(0.25) \ar@{=>} ?(0.5)_{\labelstyle \bde}
    \POS"04";"1"**{}, ?(0.25) \ar@{=>} ?(0.5)_{\labelstyle \abe}
    \or
    \POS"2" \ar "4"|{\labelstyle \ce}="24"
    \POS"0" \ar "2"|{\labelstyle \ac}="02"
    \POS"02";"1"**{}, ?(0.3) \ar@{=>} ?(0.7)^{\labelstyle \abc}
    \POS"04";"2"**{}, ?(0.2) \ar@{=>} ?(0.35)_{\labelstyle \ace}
    \POS"24";"3"**{}, ?(0.2) \ar@{=>} ?(0.6)^{\labelstyle \cde}
    \or
    \POS"1" \ar "3"|{\labelstyle \bd}="13"
    \POS"0" \ar "3"|{\labelstyle \ad}="03"
    \POS"04";"3"**{}, ?(0.2) \ar@{=>} ?(0.4)_{\labelstyle \ade}
    \POS"13";"2"**{}, ?(0.3) \ar@{=>} ?(0.7)_{\labelstyle \bcd}
    \POS"03";"1"**{}, ?(0.25) \ar@{=>} ?(0.5)^{\labelstyle \abd}
    \or
    \POS"2" \ar "4"|{\labelstyle \ce}="24"
    \POS"1" \ar "4"|{\labelstyle \be}="14"
    \POS"24";"3"**{}, ?(0.2) \ar@{=>} ?(0.6)^{\labelstyle \cde}
    \POS"04";"1"**{}, ?(0.25) \ar@{=>} ?(0.5)_{\labelstyle \abe}
    \POS"14";"2"**{}, ?(0.25) \ar@{=>} ?(0.5)^{\labelstyle \bce}
    \else\fi
  }
}

\newcommand{\pentofpent}[1]{
  \def\baselen{#1}
  \begin{xy}
    0;<\baselen,0mm>:
    *{\xybox{
        \POS(0,-4)*[o]{\pent 0}="zero"
        \POS(16,40)*[o]{\pent 3}="three"
        \POS(72,40)*[o]{\pent 1}="one"
        \POS(88,-4)*[o]{\pent 4}="four"
        \POS(44,-36)*[o]{\pent 2}="two"
        \ar@<1ex>"zero";"three"^-{\objectstyle\abcd}
        \ar@<1ex>"three";"one"^-{\objectstyle\abde}
        \ar@<1ex>"one";"four"^-{\objectstyle\bcde}
        \ar@<-1ex>"zero";"two"_-{\objectstyle\acde}
        \ar@<-1ex>"two";"four"_-{\objectstyle\abce}
        \ar@{=>}(44,-5);(44,+15)^{\objectstyle\abcde}
     }}
  \end{xy}
}

\setlist{}
\setenumerate{leftmargin=*,labelindent=\parindent}
\setitemize{leftmargin=*,labelindent=0.5\parindent}
\setdescription{leftmargin=1em}

\swapnumbers

\theoremstyle{plain}

\newtheorem{thm}{Theorem}[subsection]
\newtheorem{lem}[thm]{Lemma}
\newtheorem{cor}[thm]{Corollary}

\newtheorem{prop}[thm]{Proposition}

\theoremstyle{definition}
\newtheorem{defn}[thm]{Definition}
\newtheorem{ex}[thm]{Example}

\theoremstyle{remark}
\newtheorem{obs}[thm]{Observation}
\newtheorem{rec}[thm]{Recall}
\newtheorem{rmk}[thm]{Remark}

\makeatletter
\let\c@equation\c@thm
\makeatother
\numberwithin{equation}{subsection}

\begin{asydef}
  import "../common/adjunction" as adjunction;
\end{asydef}

\newcommand{\refI}[1]{I.\ref*{found:#1}}
\newcommand{\refII}[1]{II.\ref*{cohadj:#1}}

\title{Working Title}

\author[Riehl]{Emily Riehl}
\address{
  Department of Mathematics \\
  Harvard University \\
  Cambridge, MA 02138\\
  USA
}
\email{eriehl@math.harvard.edu}

\author[Verity]{Dominic Verity}
\address{
  Centre of Australian Category Theory \\
  Macquarie University \\
  NSW 2109 \\
  Australia
}
\email{dominic.verity@mq.edu.au}

\date{\today}

\subjclass[2010]{%
  Primary  18G55, 55U35, 55U40; %
  Secondary 18A05, 18D20, 18G30, 55U10%18D35, 18F99
}

\expandafter\ifx\csname preamble\endcsname\relax\else
  \title[Completeness results for quasi-categories of algebras]{Completeness results for quasi-categories of algebras, homotopy limits, and related general constructions}
  \date{\today}
  \renewcommand{\thethm}{\arabic{section}.\arabic{thm}}
  \numberwithin{equation}{section}
  \numberwithin{thm}{section}
  
\fi

\begin{document}

  \ifpdf
  \DeclareGraphicsExtensions{.pdf, .jpg, .tif}
  \else
  \DeclareGraphicsExtensions{.eps, .jpg}
  \fi

  \expandafter\ifx\csname abstracttext\endcsname\relax\else
  \begin{abstract}
    \abstracttext
  \end{abstract}
  \fi

  \maketitle

\tableofcontents

\fi

\section{Introduction}

Let $\qCat_\infty$ denote the simplicial category of quasi-categories, a full subcategory of the cartesian closed category of simplicial sets. Because the quasi-categories are the fibrant objects in a model category enriched over the Joyal model structure (in this case over itself), the simplicial category $\qCat_\infty$ is closed under a certain class of weighted limits: those whose weights are projective cofibrant simplicial functors. Here, a subcategory of an enriched category is \emph{closed} under a certain class of limits if it possesses those limits and if moreover they are preserved by the subcategory inclusion. 

Our aim in this paper is to prove the following theorem:

\begin{thm}\label{thm:limits-in-limits} Let $X$ be a simplicial set.  The quasi-categorically enriched subcategory of $\qCat_\infty$ spanned by those quasi-categories admitting (co)limits of shape $X$ and those functors preserving them is closed in $\qCat_\infty$ under all projective cofibrant weighted limits.
\end{thm}

This paper is a continuation of  \cite{RiehlVerity:2012tt} and \cite{RiehlVerity:2012hc}, which contain all the necessary preliminaries; references therein will have the form I.x.x.x. or II.x.x.x. Projective cofibrant weights are discussed in \S\refII{sec:weighted}, and a recognition theorem is proven. This class of weights plays a fundamental role in our development of the quasi-categorical monadicity theorem, which constructs the quasi-category of algebras for a homotopy coherent monad and characterises it up to equivalence. 

The class of \emph{projective cofibrant weighted limits}, that is, weighted limits with projective cofibrant weights, includes familiar Bousfield-Kan-style homotopy limits such as comma quasi-categories or mapping cocylinders (see Example \refII{ex:homotopyweighted}).  A weight is projective cofibrant if and only if it can be built as a retract of a countable composite of pushouts of coproducts of basic ``cells'' defined by tensoring representable functors with the boundary inclusions $\boundary\Del^n\inc\Del^n$ of simplicial sets. To prove that projective cofibrant weighted limits of quasi-categories are again quasi-categories, it suffices to show that  $\qCat_\infty$ is closed under the corresponding limit notions---namely, splittings of idempotents, limits of towers of isofibrations, pullbacks of isofibrations, products,  and cotensors with simplicial sets---and that moreover the cotensor of a quasi-category with a monomorphism is an isofibration. Here an \emph{isofibration} is a functor between quasi-categories that is a fibration in the Joyal model structure; see \refI{rec:qmc-quasicat}. Section \ref{sec:weighted} contains a general discussion of projective cofibrant weighted limits. 

Our first key result, Theorem~\ref{thm:limits-in-limits}, shows that if each quasi-category in a simplicial diagram has (co)limits of a fixed shape and each functor in that diagram preserves these, then the projective weighted limit again posseses these (co)limits and they are preserved by the functors defining the limit cone. For example, the quasi-category of algebras $A[t]$ for a homotopy coherent monad on a quasi-category $A$ is defined in %HERE: previously cited sec:monadicity, which is wrong
\S\refII{sec:formal} via a projective cofibrant weighted limit. Consequently, as an immediate corollary of our first key result, we can show that the forgetful functor $u^t\colon A[t] \to A$ creates any colimits that are preserved by the functor part $t \colon A \to A$ of the monad---precisely as is the case in ordinary category theory. %HERE: new sentence; separated into two paragraphs; first sentence in next replaces "that result" by "Theorem 1.1"
In particular, the forgetful functor creates colimits of $u^t$-split simplicial objects, completing the proof of the quasi-categorical monadicity theorem described in \S\refII{sec:monadicity} (see Corollary \ref{cor:monadic-creates-u-split}).

However, the immediate application of Theorem \ref{thm:limits-in-limits} to the question of the creation of limits by $u^t\colon A[t]\to A$ does not yield a true generalisation of the corresponding classical result, since that theorem does not require any assumption that the functor part of our monad should preserve limits. Our second key result,  %HERE: used to point to thm:monadic-limit
 Theorem~\ref{thm:monadic-completeness}, rectifies this deficiency. It makes a much more detailed and concrete analysis of the specific weight used in building quasi-categories of algebras to show that $u^t\colon A[t]\to A$ does indeed create any limits that exist in $A$ without having to constrain our homotopy coherent monad by any continuity assumptions.

In 2-category theory, the classes of {\em flexible\/} and {\em PIE\/} limits (cf.~\cite{Lack:2010fk}) are of great importance. They are, for example, the 2-limits that are inherited by 2-categories of algebras and pseudo-morphisms for a 2-monad and are therefore the ones that exist in common 2-algebraic situations such as in 2-categories of (co)complete categories, monoidal categories, Grothendieck toposes, accessible categories, and so forth. They also enjoy the familiar homotopy theoretic property that pseudo-transformations of diagrams whose legs are equivalences induce equivalences of such limits. These classes are not only parsimonious but they also encompass most 2-limit types of practical importance, such as pseudo, lax and oplax limits, comma categories, cotensors, gluing constructions, descent and lax descent constructions, and so forth. The class of limits weighted by projective cofibrant functors or projective cell complexes are direct quasi-categorical analogues of the classes of flexible and PIE limits respectively. They too contain analogues of all of the 2-limit types discussed above and much more besides. Consequently, the completeness results we develop in this paper are ripe for wide application.

The techniques we use to prove Theorem \ref{thm:limits-in-limits} tell us how the limits or colimits in a projective cofibrant weighted limit are constructed, enabling more precise applications. For example, a quasi-category is \emph{stable} if it is finitely complete and cocomplete and if moreover the terminal object is initial and a square is a pullback if and only if it is a pushout \cite[1.1.3.4]{Lurie:2012uq}. A functor of stable quasi-categories is \emph{exact} if it preserves finite limits or finite colimits, either condition implying the other \cite[1.1.4.1]{Lurie:2012uq}. Theorem \ref{thm:limits-in-limits} immediately implies that the homotopy limit of a diagram of stable quasi-categories and exact functors is finitely complete and cocomplete, but we can show further that the homotopy limit is again a stable quasi-category; see Theorem \ref{thm:stable}.

Our interest in the questions discussed here originally arose from our work on homotopy coherent monads presented in~\cite{RiehlVerity:2012hc}. We had also been interested in some specific (co)completeness results for functor quasi-categories (hom-spaces in $\qCat_\infty$), comma constructions, homotopy limits, and the like, but had not sought to unify the results we had developed. However, an email sent to us by Tom Fiore encouraged us to look at the more general results we present here. Specifically, he asked whether it was possible to select certain families of limits in a quasi-category in a suitably functorial manner. Our approach to this problem was to observe that his result could be reduced to asking whether a particular limit of quasi-categories weighted by  some projective cofibrant weight possessed limits of that kind. The results presented here show that it does.

\subsection*{Limits and colimits in a quasi-category}

In \S\refI{sec:limits}, we show that limits (resp.~colimits) in a quasi-category are characterised by certain absolute right (resp.~left) lifting diagrams in $\qCat_2$, Andr\'{e} Joyal's 2-category of quasi-categories. Absolute right liftings are introduced in Definition \refI{defn:abs-right-lift} and recalled in \S\ref{sec:absolute}, where precise details are first needed. For ease of exposition, we restrict our discussion to limits from here on; the corresponding results for colimits are dual. 

{
\renewcommand{\thethm}{\refI{defn:families.of.diagrams}}
\begin{defn}
A quasi-category $A$ \emph{has all limits of shape} $X$ if and only if there is an absolute right lifting  
  \begin{equation*}
      \xymatrix{ \ar@{}[dr]|(.7){\Downarrow\lambda} & A \ar[d]^c \\ A^X \ar[ur]^{\lim} \ar[r]_{\id_{A^X}} & A^X}
  \end{equation*} 
  of the identity through the constant diagram functor $c \colon A \to A^X$. A functor $f \colon A \to B$ \emph{preserves all limits of shape} $X$ if the composite diagram 
  \begin{equation*}
      \xymatrix{ \ar@{}[dr]|(.7){\Downarrow\lambda} & A \ar[d]^c \ar[r]^f & B \ar[d]^c \\ A^X \ar[ur]^{\lim} \ar[r]_{\id_{A^X}} & A^X \ar[r]_{f^X} & B^X}
  \end{equation*} 
is again an absolute right lifting; see, for example, \refI{prop:RAPL}.
\end{defn}
\addtocounter{thm}{-1}
}

Taking $X = \emptyset$, Definition \refI{defn:families.of.diagrams} specialises to say that a quasi-category $A$ has a terminal object if and only if there is an absolute right lifting 
  \begin{equation*}
      \xymatrix{ \ar@{}[dr]|(.7){\Downarrow\lambda} & A \ar[d]^c \\ \catone \ar[ur]^{a} \ar[r]_{\id_{\catone}} & \catone}
  \end{equation*} 
  Here the 2-cell $\lambda$ is necessarily an identity because the quasi-category $\catone$ is terminal. By a standard exercise in formal category theory, this is equivalent to asking that  the constant functor $A \to \catone$ has a right adjoint $\catone \xrightarrow{a} A$ that picks out the terminal object. The unit of this adjunction is represented by a 1-simplex in $A^A$ from the identity functor to the constant functor at the terminal object whose component at the terminal vertex is degenerate; see Lemma \refI{lem:min-term-pres}.

\subsection*{Outline of the proof}

The special case of terminal objects will be of particular interest: the general case for limits of shape $X$ may be reduced to this special one.  We shall prove Theorem~\ref{thm:limits-in-limits} in two steps, which occupy sections \ref{sec:terminal} and \ref{sec:absolute} respectively:
  \begin{enumerate}[label=(\Roman*)]
    \item\label{item:stepI} Prove the special case of  quasi-categories which admit terminal objects.
    \item\label{item:stepII} Reduce the general problem to that special case ``fibre-wise.'' 
  \end{enumerate}
  
The key ingredient for this second step is Theorem \refI{thm:pointwise}, recalled in \S\ref{sec:absolute}, which encodes absolute right lifting diagrams as a collection of terminal objects in particular comma quasi-categories. In the special case of absolute lifting diagrams of the form given in Definition \refI{defn:families.of.diagrams}, the comma quasi-categories of Theorem \refI{thm:pointwise} are quasi-categories of cones over a fixed diagram. A terminal object encodes a limiting cone; see Observation \refI{obs:limits.are.limits}.

Special cases of Theorem \ref{thm:limits-in-limits} have been established before: for example, \cite[5.4.5.5]{Lurie:2009fk} proves this result for homotopy pullbacks of quasi-categories admitting and functors preserving particular colimits; see also \cite[5.5.7.11]{Lurie:2009fk}. But we know of no other treatment at this level of generality. Incidentally, the strategy used in \cite[\S 5.4.5]{Lurie:2009fk} is similar---the case of general colimits is deduced from the special case of initial objects---although the details of the proof are quite different. 

\subsection*{Acknowledgments} 

We would like to thank Tom Fiore for encouraging us to complete the work presented here. This material is based upon work supported by the National Science Foundation under Award No.~DMS-1103790 and by the Australian Research Council under Discovery grant number DP130101969.
  
\section{Weighted limits in \texorpdfstring{$\protect\qCat_\infty$}{qCat}}\label{sec:weighted}

The theory of weighted limits is reviewed in \S\refII{subsec:weighted}.  For the reader's convenience, and as a warm-up to our proof of Theorem \ref{thm:limits-in-limits}, recall:
{
\renewcommand{\thethm}{\refII{prop:projwlims}}
\begin{prop}
The full simplicial subcategory $\qCat_\infty$ of quasi-categories is closed  in $\sSet$ under each of the following classes of limits:
 \begin{enumerate}[label=\arabic*)]
    \item arbitrary small products,
    \item pullbacks of isofibrations,
    \item countable composites of isofibrations, 
    \item splitting of idempotents, and
    \item cotensors with any simplicial set.
\end{enumerate}
Moreover the cotensor of a quasi-category with a monomorphism of simplicial sets is an isofibration.  Hence, $\qCat_\infty$ is closed under limits weighted by projective cofibrant weights.
\end{prop}
\addtocounter{thm}{-1}
}
\begin{proof}
Let $\stcat{A}$ be a small simplicial category. A weight $W \in \sSet^{\stcat{A}}$ is projective cofibrant just when it is expressible as a retract of a countable composite of pushouts of coproducts of maps $i \times \stcat{A}_a \colon X \times \stcat{A}_a \to Y \times \stcat{A}_a$ built by tensoring a covariant representable $\stcat{A}_a$ with an inclusion $X\inc Y$ of simplicial sets. By cocontinuity of the weighted limit bifunctor in the weight, the weighted limit $\wlim{W}{D}_{\stcat{A}}$ of a diagram $D \colon \stcat{A} \to \qCat_\infty$ is a retract of a countable tower of pullbacks of products of maps \[ \wlim{i \times \stcat{A}_a}{D}_{\stcat{A}} \cong Da^i \colon Da^Y \to Da^X.\] 

Because the quasi-categories are the fibrant objects in a monoidal model structure with all objects cofibrant, $\qCat_\infty$ is closed under simplicial cotensors, and in particular the domains and codomains of these maps are again quasi-categories. Moreover, monoidalness of the Joyal model structure implies that each of these maps is an isofibration (a fibration between fibrant objects).

As a category of fibrant objects, the quasi-categories are closed under products and splittings of idempotents as well as pullbacks and composites of isofibrations. In this way, we see that the weighted limit $\wlim{W}{D}_{\stcat{A}}$ is again a quasi-category.
\end{proof}

\begin{rmk}
In the proof of Proposition \refII{prop:projwlims} just given, we have used a slightly different definition of projective cofibrant weight than appears in \refII{defn:proj-cof}. These are of course equivalent: by the algebraic small object argument, a cofibrant object in the projective model structure on $\sSet^\stcat{A}$ is expressible as a retract of a countable composite of pushouts of coproducts of maps of the form $\partial\Del^n \times \stcat{A}_a \to \Del^n\times\stcat{A}_a$. Only countable composites are needed because the simplicial sets appearing here are $\omega$-small.
\end{rmk}

\begin{ex}
The class of projective cofibrant weighted limits includes:
\begin{itemize}
\item cotensors by arbitrary simplicial sets (\refII{ex:diagrams-weighted-limit}),
\item comma objects and iso-comma objects (\refII{ex:commaweightedlimit}),
\item the appropriate quasi-categorical versions of inserters and equifiers of 2-category theory (left as an exercise for the interested reader),
\item Bousfield-Kan-style homotopy limit for diagrams of any shape (\refII{ex:homotopyweighted}),
\item algebras for a homotopy coherent monad (\refII{lem:projcof1}; see also \ref{rec:algebras} below).
\end{itemize}
\end{ex}

\section{Quasi-categories with terminal objects}\label{sec:terminal}

The basic theory of terminal objects is developed in \S\refI{subsec:terminal}. We begin by introducing some special notation for the quasi-categorically enriched category described in the statement of Theorem~\ref{thm:limits-in-limits}.

\begin{defn}\label{defn:qCat-X}
  For any simplicial set $X$, let $\qCat_{X\infty}$ denote the quasi-categorically enriched subcategory of $\qCat_\infty$ with:
  \begin{itemize}
    \item {\bf objects} the quasi-categories $A$ which admit all limits of shape $X$,
    \item {\bf $n$-arrows} all those $n$-arrows whose vertices are functors $f\colon A\to B$ which preserve all limits of shape $X$.
  \end{itemize}
  We shall write $\qCat_X$ for the underlying category of $\qCat_{X\infty}$.
\end{defn}

Our aim in this section is to prove that $\qCat_{\emptyset\infty}$, the subcategory of quasi-categories admitting terminal objects and terminal-object-preserving functors, is closed in $\qCat_\infty$ under limits with projective cofibrant weights.

\begin{obs}\label{obs:terminal-pres-banal}
In order to show that a functor $f\colon A\to B$ of quasi-categories preserves terminal objects, it is enough to show that there is {\em some\/} terminal object $t$ in $A$ such that $ft$ is terminal in $B$: all terminal objects in $A$ are isomorphic, and any functor between quasi-categories preserves isomorphisms.
\end{obs}

\begin{obs}\label{obs:terminals-in-limits-prf}
  Our first order of work will be to prove that $\qCat_\emptyset$ is closed in $\qCat$ under each of the following classes of conical limits:
  \begin{enumerate}[label=\arabic*)]
    \item arbitrary small products,
    \item pullbacks of isofibrations,
    \item countable composites of isofibrations, and
    \item splitting of idempotents.
  \end{enumerate}
The ``isofibrations'' referred to here are isofibrations that preserve terminal objects, i.e., isofibrations that lie in the subcategory $\qCat_\emptyset$. We shall prove these closure properties in a series of lemmas. In Corollary \ref{cor:terminal-cotensor}, we will  prove that the quasi-categorically enriched category $\qCat_{\emptyset\infty}$ is closed in $\qCat_\infty$ under:
  \begin{enumerate}[resume*]
    \item cotensors with any simplicial set.
  \end{enumerate}
Moreover, we will show that the map induced by cotensoring with an inclusion of simplicial sets is an isofibration in $\qCat_\emptyset$. As in the proof of Proposition \refII{prop:projwlims}, it will follow that  $\qCat_{\emptyset\infty}$ is closed in $\qCat_\infty$ under limits weighted by projective cofibrant weights.
\end{obs}

\begin{rmk}[proof-schema]\label{rmk:qCat-term-closure-expl}
  For each class of conical limits listed in Observation~\ref{obs:terminals-in-limits-prf}, we shall prove the desired closure result using a similar argument. Specifically, given a diagram in $\qCat_\emptyset$ we shall take its limit $L$ in $\qCat$, which comes equipped with projection functors $\pi_i\colon L\to A_i$, and establish the following results:
  \begin{enumerate}[label=(\roman*)]
    \item\label{enum:term-in-lim-p1} there is some object $t$ of $L$ which is {\em pointwise terminal}, in the sense that each $\pi_it$ is terminal in $A_i$, and
    \item\label{enum:term-in-lim-p2} every pointwise terminal object in $L$ is a terminal object in $L$.
  \end{enumerate} 
Now~\ref{enum:term-in-lim-p1} delivers a pointwise terminal $t$ which is terminal in $L$ by~\ref{enum:term-in-lim-p2}. The pointwise terminality of $t$ tells us that each $\pi_i$ maps that particular terminal to a terminal in $A_i$, so it follows that $\pi_i$ preserves all terminal objects by Observation~\ref{obs:terminal-pres-banal}. Together these arguments show that the vertex $L$ and its limiting cone are in $\qCat_\emptyset$. 

  Finally, if $f\colon A\to L$ is a functor induced by another cone in $\qCat_\emptyset$ then we know that the components $f_i=\pi_if\colon A\to A_i$ are all in $\qCat_\emptyset$ and thus preserve terminal objects. It follows that a terminal object $s$ of $A$ maps to an object $fs$ of $L$ for which $\pi_ifs = f_is$ is terminal in $A_i$ for each $i$. In other words, $fs$ is pointwise terminal in $L$ and hence is thus terminal by~\ref{enum:term-in-lim-p2}. This demonstrates  that $f$ preserves terminal objects as required.
\end{rmk}

\subsection*{Universal properties of terminal objects}

We will make repeated use of a few simple lemmas expressing the (relative) universal property of terminal objects and their behaviour with respect to isofibrations.

\begin{rec}[terminal objects]\label{rec:terminal.char}
  Propositions~\refI{prop:terminalconverse} and~\refI{prop:terminal-ext-univ} tell us that if $a$ is an object of a quasi-category $A$ then the following are equivalent:
  \begin{itemize}
    \item The vertex $a$ is a terminal object of $A$ in the sense of Joyal~\cite{Joyal:2002:QuasiCategories}:  it enjoys the universal property that any sphere $\boundary\Del^n\to A$, $n\geq 1$, which maps the final vertex $\fbv{n}$ to $a$ 
    \begin{equation*}
\xymatrix{ \Del^0 \ar[r]_-{\fbv{n}} \ar@/^1pc/[rr]^a & \boundary\Del^n \ar[r] \ar@{u(->}[d] & A  \\  & \Del^n \ar@{-->}[ur] & }\end{equation*}
    has a filler, this being the dashed arrow in the diagram above.
    \item For all quasi-categories $B$, the composite functor $\xymatrix@1{{B}\ar[r]^{!} & {\Del^0}\ar[r]^a & {A}}$ is terminal in the hom-category $\hom_2(B,A)=\ho(A^B)$ of $\qCat_2$.
    \item The functor $a\colon\Del^0\to A$ is right adjoint to the unique functor $!\colon A\to \Del^0$.
  \end{itemize}
\end{rec}

{
\renewcommand{\thethm}{\refII{lem:rel-lift-terminal}}
\begin{lem}
  Suppose that $E$ and $B$ are quasi-categories which possess terminal objects and that $p\colon E\tfib B$ is an isofibration which preserves terminal objects. Assume also that $e$ is terminal in $E$. Then any lifting problem
  \begin{equation*}
    \xymatrix@=2em{
      \Del^0 \ar[r]_-{\fbv{n}} \ar@/^1pc/[rr]^e &
      {\boundary\Del^n}\ar[r]_-u\ar@{u(->}[d] &
      {E}\ar@{->>}[d]^p \\
       & {\Del^n}\ar@{-->}[ur]\ar[r]_-v & {B}
    } 
  \end{equation*}
  with $n>0$ in which $u$ carries the vertex $\fbv{n}$ to $e$ has a solution.
\end{lem}
\addtocounter{thm}{-1}
}

\begin{obs}
  We observe in passing that the lifting property of the last lemma is equivalent to the statement that 
  \begin{equation*}
    \xymatrix@=1.5em{
      {}\ar@{}[dr]|(0.68){||} & {E}\ar@{->>}[d]^p \\
      {\catone}\ar[ur]^{e}\ar[r]_{pe} & {B}
    }
  \end{equation*}
  is an absolute right lifting diagram in $\qCat_2$.
\end{obs}

Lemma \refII{lem:rel-lift-terminal} possesses the following immediate converse:

\begin{lem}\label{lem:rel-lift-terminal-conv}
  Suppose that $p\colon E\tfib B$ is an isofibration of quasi-categories and that $e$ is an object of $E$ for which $pe$ is a terminal object in $B$. Assume also that the lifting property of Lemma~\refII{lem:rel-lift-terminal} holds. Then $e$ is a terminal object of $E$. 
\end{lem}

\begin{proof}
  Suppose that $u\colon\boundary\Del^n\to E$ is a sphere in $E$ which carries the vertex $\fbv{n}$ to the object $e$. Then the sphere $pu\colon\boundary\Del^n\to B$ carries $\fbv{n}$ to the terminal object $pe$ of $B$ and so it has a filler $v\colon\Del^n\to B$. This provides us with a lifting problem of the form given in the statement of Lemma~\refII{lem:rel-lift-terminal}, which we may solve to provide a filler in $E$ for the original sphere $u$,  establishing the universal property which shows that $e$ is terminal in $E$.
\end{proof}

\begin{lem}\label{lem:rel-lift-terminal-companion}
  Suppose that $E$ and $B$ are quasi-categories which possess terminal objects and that $p\colon E\tfib B$ is an isofibration which preserves terminal objects. If $s$ is a terminal object in $B$ then there exists a terminal object $t$ in $E$ with $pt = s$.
\end{lem}

\begin{proof}
The hypotheses imply that $E$ admits a terminal object $t'$ and that $pt'$ is terminal in $B$. By uniqueness of terminal objects, there exists some isomorphism $s\cong pt'$ in $B$ which we can lift along the isofibration $p$ to give $t\cong t'$ which lies over $s\cong pt'$. It follows that $t$ is also terminal in $E$, because it is isomorphic to our original terminal object $t'$, and that $pt = s$ as required.
\end{proof}

\subsection*{Conical limits in \texorpdfstring{$\protect\qCat_\emptyset$}{qCat}}

We will now show that $\qCat_\emptyset$ is closed under the conical limits of Observation \ref{obs:terminals-in-limits-prf}.

\begin{lem}\label{lem:qCat_t-closed-products}
    The subcategory $\qCat_\emptyset$ is closed in $\qCat$ under arbitrary small products.
\end{lem}

\begin{proof} 
  Let $\{A_i\}_{i\in I}$ be an indexed family of quasi-categories in $\qCat_\emptyset$, and let $\prod_{i\in I} A_i$ denote their product in $\qCat$ with projection functors $\pi_i\colon \prod_{i\in I} A_i \to A_i$. We establish the properties discussed in Remark \ref{rmk:qCat-term-closure-expl}:
  \begin{enumerate}[label=(\roman*)]
    \item Any choices of terminal object $t_i$ in $A_i$ for each $i\in I$ can be collected together to give a pointwise terminal $(t_i)_{i\in I}$ in $\prod_{i\in I} A_i$.
    \item Suppose that $(t_i)_{i\in I}$ is pointwise terminal and that we have a lifting problem
    \begin{equation}\label{eq:to-solve-prod}
      \xymatrix@R=1.5em@C=2.5em{
        {\boundary\Del^n}\ar[r]^-{f}\ar@{u(->}[d] &
        {\prod_{i\in I} A_i} \\
        {\Del^n}\ar@{-->}[ur]_-{?}
      }
    \end{equation}
    in which $f$ maps $\fbv{n}$ to $(t_i)_{i\in I}$. Composing \eqref{eq:to-solve-prod} with the projection $\pi_i$ we get a lifting problem
    \begin{equation*}
      \xymatrix@R=1.5em@C=2.5em{
        {\boundary\Del^n}\ar[r]^-{f}\ar@{u(->}[d] & \prod_{i \in I} A_i \ar[r]^-{\pi_i} &
        {A_i} \\
        {\Del^n}\ar@{-->}[urr]_{\bar{f}_i}
      }
    \end{equation*}
    in which $\pi_i f$ maps $\fbv{n}$ to the terminal $t_i$ in $A_i$. By the universal property of the terminal $t_i$, as recalled in~\ref{rec:terminal.char}, we may find a solution $\bar{f}_i$ to this problem as displayed. The collection of all of these lifts induces a functor $\bar{f}\colon \Del^n\to\prod_{i\in I} A_i$ with $\pi_i f = f_i$, which provides us with a solution to the original lifting problem~\eqref{eq:to-solve-prod}. This proves that the object $(t_i)_{i \in I}$ also possesses the universal property of \ref{rec:terminal.char}, making it terminal in $\prod_{i \in I} A_i$.\qedhere
  \end{enumerate}
\end{proof}

\begin{lem}\label{lem:qCat_t-closed-pullbacks}
    The subcategory $\qCat_\emptyset$ is closed in $\qCat$ under pullbacks of isofibrations.
\end{lem}

\begin{proof}
Given an isofibration $p\colon E\tfib B$ and a functor $f\colon A\to B$ in $\qCat_\emptyset$, we form the pullback
  \begin{equation*}\label{eq:pb-in-qCat}
    \xymatrix@=2em{
      {E\times_B A}\pbexcursion
      \ar[r]^-{\pi_0}\ar@{->>}[d]_{\pi_1} &
      {E}\ar@{->>}[d]^{p} \\
      {A}\ar[r]_{f} & {B}
    }
  \end{equation*}
  in $\qCat$. To show that $E \times_B A$ has and $\pi_1$ preserves terminal objects, we follow the outline of Remark~\ref{rmk:qCat-term-closure-expl}:
  \begin{enumerate}[label=(\roman*)]
    \item Given a terminal object $t \in A$, $ft$ is terminal in $B$, and we may apply Lemma~\ref{lem:rel-lift-terminal-companion} to obtain a terminal object $s$ in $E$ with $ps = ft$. The pair $(s,t)$ is an object of $E\times_B A$ which is pointwise terminal by construction.
    \item Suppose we are given a pointwise terminal $(s,t)$ in $E\times_B A$ and a lifting problem
    \begin{equation}\label{eq:to-solve-pb}
      \xymatrix@R=1.5em@C=2.5em{
        {\boundary\Del^n}\ar[r]^-{g}\ar@{u(->}[d] &
        {E\times_B A} \\
        {\Del^n}\ar@{-->}[ur]_-{?}
      }
    \end{equation}
    in which $g$ maps $\fbv{n}$ to $(s,t)$. To solve this problem consider the following diagram
  \begin{equation*}
    \xymatrix@=2em{
      {\boundary\Del^n}\ar[r]^-{g}\ar@{u(->}[d] &
      {E\times_B A}
      \ar[r]^-{\pi_0}\ar@{->>}[d]_(0.3){\pi_1} &
      {E}\ar@{->>}[d]^{p} \\
      {\Del^n}\ar@{-->}[r]_{\bar{g}_1}
      \ar@{-->}[urr]_(0.7){\bar{g}_0} &
      {A}\ar[r]_{f} & {B}
    }
  \end{equation*}
  in which the dashed arrow $\bar{g}_1$ is a solution to the lifting problem $\pi_1g\colon\boundary\Del^n\to A$; this solution exists because $\pi_1g$ maps $\fbv{n}$ to $t$, which is terminal in $A$.

The dashed solution $\bar{g}_0$ exists by application of Lemma~\refII{lem:rel-lift-terminal}, which applies since $p$ is a terminal-object-preserving isofibration and $\pi_0 g$ maps $\fbv{n}$ to the  terminal object $s$ of $E$. The lower right triangle in this diagram tells us that $p\bar{g}_0=f\bar{g}_1$, so the pair $(\bar{g}_0,\bar{g}_1)$ induces a map $\bar{g}\colon\Del^n\to E\times_B A$ which provides a solution to the original lifting problem~\eqref{eq:to-solve-pb} as required.\qedhere
  \end{enumerate}
\end{proof}

\begin{lem}\label{lem:qCat_t-closed-countable-composite}
  The subcategory $\qCat_\emptyset$ is closed in $\qCat$ under countable composites of isofibrations.
\end{lem}

\begin{proof}
  Suppose that we are given a countable sequence 
  \begin{equation*}
    \xymatrix{
      {}\ar@{.}[r] & {}\ar@{->>}[r]^-{p_i} & 
      {A_i} \ar@{->>}[r]^{p_{i-1}} & {} \ar@{.}[r] &
      {} \ar@{->>}[r]^-{p_3} & 
      {A_3} \ar@{->>}[r]^{p_2} & 
      {A_2} \ar@{->>}[r]^{p_1} & 
      {A_1} \ar@{->>}[r]^{p_0} & {A_0}
    }  
  \end{equation*}
  of quasi-categories and isofibrations in $\qCat_\emptyset$, and let $A$ be their limit in $\qCat$. We know that the projection functors $\pi_i\colon A\tfib A_i$ are isofibrations. Once again, we establish the conditions of Remark~\ref{rmk:qCat-term-closure-expl}:
  \begin{enumerate}[label=(\roman*)]
    \item Our task is to build a family $(t_i)_{i\in\mathbb{N}}$ of terminal objects $t_i$ in $A_i$ which satisfies the compatibility property that $p_it_{i+1} = t_i$ for all $i\in\mathbb{N}$. To do this we simply proceed inductively, picking any terminal object $t_0$ of $A_0$ and then extending our family from level $i$ to level $i+1$ by applying Lemma~\ref{lem:rel-lift-terminal-companion} to the terminal object-preserving-isofibration $p_i\colon A_{i+1}\tfib A_i$  in order to find a $t_{i+1}$ which is terminal in $A_{i+1}$ and which has $p_it_{i+1}=t_i$.
    \item Suppose $t$ in $A$ is a pointwise terminal object projecting to the family $(t_i)_{i \in\mathbb{N}}$. To show that $t$ is terminal, we must solve any lifting problem
    \begin{equation}\label{eq:to-solve-tower}
          \xymatrix@R=1.5em@C=2.5em{
        {\boundary\Del^n}\ar[r]^-{f}\ar@{u(->}[d] &
        {A} \\
        {\Del^n}\ar@{-->}[ur]_-{?}
      }
    \end{equation}
    that maps the vertex $\fbv{n}$ to $t$. We construct a solution $\bar{f}$ by inductively defining its components. Start by using the fact that $t_0$ is terminal to find a lift
        \begin{equation*}
          \xymatrix@R=1.5em@C=2.5em{
        {\boundary\Del^n}\ar[r]^-{f}\ar@{u(->}[d] &
        {A} \ar[r]^-{\pi_0} & A_0 \\
        {\Del^n}\ar@{-->}[urr]_-{\bar{f}_0}
      }
    \end{equation*}
Suppose now that we have defined maps $\bar{f}_i\colon\Del^n\to A_i$ so that $p_i\bar{f}_{i+1} = \bar{f}_i$ for all $i < k$. By  Lemma~\refII{lem:rel-lift-terminal} and the hypothesis that $t_k$ is terminal, there is a solution to the lifting problem 
        \begin{equation*}
          \xymatrix@R=1.5em@C=2.5em{
        {\boundary\Del^n}\ar[r]^-{f}\ar@{u(->}[d] &
        {A} \ar[r]^-{\pi_k} & A_k \ar@{->>}[d]^{p_{k-1}} \\
        {\Del^n}\ar@{-->}[urr]_-{\bar{f}_k} \ar[rr]_{\bar{f}_{k-1}} & & A_{k-1}
      }
    \end{equation*}
    defining $\bar{f}_k$. The family $(\bar{f}_i)_{i \in\mathbb{N}}$ defines the desired solution $\bar{f}$ to \eqref{eq:to-solve-tower}.\qedhere
  \end{enumerate}
\end{proof}

\begin{lem}\label{lem:qCat_t-closed-idempotents}
  The subcategory $\qCat_\emptyset$ is closed in $\qCat$ under splitting of idempotents.
\end{lem}

\begin{proof}
  Suppose that $e\colon A\to A$ is an idempotent in $\qCat_\emptyset$ and let $A^e$ denote the maximal sub-quasi-category of $A$ stabilised by that idempotent---i.e., the equaliser of $e$ and $\id_A$---whose single limit cone projection is the inclusion $A^e\inc A$. We verify the conditions of Remark~\ref{rmk:qCat-term-closure-expl} as follows:
  \begin{enumerate}[label=(\roman*)]
  \item For any terminal object $t$ of $A$, the object $et$ of $A^e$ is terminal in $A$, since $e$ is assumed to preserve terminal objects. Thus $et$ is pointwise terminal in $A^e$.
  \item Now suppose that $g\colon\boundary\Del^n\to A^e$ is a lifting problem which maps $\fbv{n}$ to an object $s$ which is pointwise terminal in $A^e$. This means that $s$ is terminal in $A$, which ensures the existence of the dashed solution in the following lifting problem:
  \begin{equation*}
    \xymatrix@R=1.5em@C=2.5em{
      {\boundary\Del^n}\ar[r]^-{g}\ar@{u(->}[d] &
      {A^e}\ar@{u(->}[d] \\
      {\Del^n}\ar@{-->}[r]_-{\bar{g}} & {A}
    }
  \end{equation*}
  It is now easy to check that $e\bar{g}\colon\Del^n\to A^e$ is a solution to our original lifting problem defined by the sphere $g\colon\boundary\Del^n\to A^e$.\qedhere
  \end{enumerate}
\end{proof}

We defer the proof that $\qCat_\emptyset$ is closed under cotensors by simplicial sets and that cotensors with monomorphisms induce terminal-object-preserving isofibrations to Lemma \ref{lem:cotensor-closure}, which proves a more general version of this statement. Modulo this step, which appears as Corollary \ref{cor:terminal-cotensor}, we may now prove our desired result.

\begin{thm}\label{thm:terminal-closure} $\qCat_{\emptyset,\infty}$ is closed in $\qCat_\infty$ under limits weighted by projective cofibrant weights.
\end{thm}
\begin{proof}
A projective cofibrant weighted limit may be expressed as a retract of the limit of a countable tower of pullbacks of products of maps defined by cotensoring an object in the diagram with a monomorphism of simplicial sets. Corollary \ref{cor:terminal-cotensor} shows that these maps are terminal-object preserving isofibrations between quasi-categories admitting terminal objects. Lemmas \ref{lem:qCat_t-closed-products}, \ref{lem:qCat_t-closed-pullbacks}, \ref{lem:qCat_t-closed-countable-composite}, and \ref{lem:qCat_t-closed-idempotents} show that the rest of the diagram lies in $\qCat_\emptyset$ as well.
\end{proof}

\section{Absolute lifting diagrams via terminal objects}\label{sec:absolute}

We now turn our attention to step \ref{item:stepII}: the reduction of the general case of Theorem \ref{thm:limits-in-limits} to the special one of quasi-categories admitting and functors preserving terminal objects. First we must set up the precise result that we intend to prove.

\subsection*{Absolute right lifting diagrams and right exact transformations}

\begin{rmk}
  As in \refI{ntn:pb.po.joins}, we let $\pbshape$ denote the category $\xymatrix@1{{c}\ar[r] & {a} & {b}\ar[l]}$ and write $\qCat_\infty^\pbshape$ for the simplicially enriched category of simplicial functors $\pbshape\to\qCat_\infty$. Objects of $\qCat_\infty^\pbshape$ are pairs of maps
  \begin{equation}\label{eq:obj-qCat-pbshape}
    \xymatrix@=1.5em{
      {} & {B}\ar[d]^f \\
      {C}\ar[r]_g & {A}
    }
  \end{equation}
  and 0-arrows are natural transformations
  \begin{equation}\label{eq:obj-qCat-pbshape-trans}
    \xymatrix@=1.5em{
      {} & {B}\ar[d]_{f}\ar[dr]^{v} & {} \\
      {C}\ar[r]^{g}\ar[dr]_{w} & {A}\ar[dr]|*+<1pt>{\scriptstyle u} & 
      {B'}\ar[d]^{f'} \\
      {} & {C'}\ar[r]_{g'} & {A'}
    }
  \end{equation}
  We say that a transformation of this kind is a {\em pointwise isofibration\/} if each of its components $u$, $v$, and $w$ is an isofibration.

  As with any enriched functor category, $\qCat_\infty^\pbshape$ inherits limits pointwise from $\qCat_\infty$. In particular, it  has limits weighted by projective cofibrant weights, which are then preserved by the simplicial projection functors $P_a,P_b,P_c\colon\qCat_\infty^\pbshape\to\qCat_\infty$  which evaluate at each of the objects $a$, $b$, and $c$ of $\pbshape$ respectively.
\end{rmk}

  Our particular interest is in a subcategory of $\qCat_\infty^\pbshape$ whose objects are those diagrams~\eqref{eq:obj-qCat-pbshape} in which $g$ admits an absolute right lifting through $f$ in Joyal's 2-category of quasi-categories $\qCat_2$.

  {
  \renewcommand{\thethm}{\refI{defn:abs-right-lift}}
  \begin{defn} 
    \addtocounter{thm}{-1}
  An \emph{absolute right lifting diagram} in a 2-category  consists of the data
    \begin{equation}\label{eq:lifting.triangle}
    \xymatrix@=1.5em{ \ar@{}[dr]|(.7){\Downarrow\lambda} & B \ar[d]^f \\ C \ar[r]_g \ar[ur]^\ell & A}
  \end{equation}
     with the universal property that if we are given any 2-cell $\chi$ of the form depicted to the left of the following equality
\[
    \vcenter{\xymatrix{ X \ar[d]_c \ar[r]^b \ar@{}[dr]|{\Downarrow\chi} & B \ar[d]^f \\ C \ar[r]_g & A}} \mkern20mu = \mkern20mu \vcenter{\xymatrix{ X \ar[d]_c \ar[r]^b \ar@{}[dr]|(.3){\exists !\Downarrow}|(.7){\Downarrow\lambda} & B \ar[d]^f \\ C \ar[ur]|(.4)*+<2pt>{\scriptstyle\ell} \ar[r]_g & A}}
 \]
 then it admits a unique factorisation of the form displayed to the right of that equality.
\end{defn}
}

\begin{defn}\label{defn:abs-lift-trans-right-exactness}
Transformations of the kind depicted in~\eqref{eq:obj-qCat-pbshape-trans} between diagrams which admit absolute right liftings give rise to the following diagram
  \begin{equation}\label{eq:induced-mate}
    \vcenter{
      \xymatrix@=1.5em{
        {} \ar@{}[dr]|(.7){\Downarrow\lambda}
        & {B}\ar[d]^{f}\ar[dr]^{v} & {} \\
        {C}\ar[ur]^{\ell}\ar[r]_{g}\ar[dr]_{w} & 
        {A}\ar[dr]|*+<1pt>{\scriptstyle u} & 
        {B'}\ar[d]^{f'} \\
        {} & {C'}\ar[r]_{g'} & {A'}
      }
    } 
    \mkern40mu = \mkern40mu
    \vcenter{
      \xymatrix@=1.5em{
        {} & {B}\ar[dr]^{v} & {} \\
        {C}\ar[ur]^{\ell}\ar[dr]_{w} & 
        {\scriptstyle\Downarrow\tau}\ar@{}[dr]|(.7){\Downarrow\lambda'} & 
        {B'}\ar[d]^{f'} \\
        {} & {C'}\ar[ur]^{\ell'}\ar[r]_{g'} & {A'}
      }
    } 
  \end{equation}
  in which the triangles are absolute right liftings and the 2-cell $\tau$ is induced by the universal property of the triangle on the right. We say that the transformation \eqref{eq:induced-mate} is {\em right exact\/} if and only if the induced 2-cell $\tau$ is an isomorphism. This right exactness condition holds if and only if, in the diagram on the left, the whiskered 2-cell $u\lambda$ displays $v\ell$ as the absolute right lifting of $g'w$ through $f'$.
\end{defn}

\begin{defn}
  Let $\qCat_{r\infty}^\pbshape$ denote the simplicial subcategory of $\qCat_\infty^\pbshape$ with:
  \begin{itemize}
    \item {\bf objects} those diagrams~\eqref{eq:obj-qCat-pbshape} in which $g$ admits an absolute right lifting through $f$, and
    \item {\bf $n$-arrows} those $n$-arrows of $\qCat_\infty^\pbshape$ whose vertices are right exact transformations.
  \end{itemize}
As usual, we let $\qCat_r^\pbshape$ denote the underlying category of $\qCat_{r\infty}^\pbshape$.
\end{defn}

\begin{obs}\label{obs:comparison-2cell-composite}
    To check that this definition makes sense, we need  to show that the composite of a pair of right exact transformations is itself right exact. This follows from the following calculation 
  \begin{equation*}
    \vcenter{
      \xymatrix@=1.5em{
        {} \ar@{}[dr]|(.7){\Downarrow\lambda}
        & {B}\ar[d]^{f}\ar[dr]^{v} & {} & {} \\
        {C}\ar[ur]^{\ell}\ar[r]_{g}\ar[dr]_{w} & 
        {A}\ar[dr]|*+<1pt>{\scriptstyle u} & {B'}\ar[d]^{f'}\ar[dr]^{v'} & {} \\
        {} & {C'}\ar[r]_{g'}\ar[dr]_{w'} & 
        {A'}\ar[dr]|*+<1pt>{\scriptstyle u'} & {B''}\ar[d]^{f''} \\
        {} & {} & {C''}\ar[r]_{g''} & {A''}
      }
    } 
    \mkern5mu = \mkern5mu
    \vcenter{
      \xymatrix@=1.5em{
        {} & {B}\ar[dr]^{v} & {} & {} \\
        {C}\ar[ur]^{\ell}\ar[dr]_{w} & 
        {\scriptstyle\Downarrow\tau}
        \ar@{}[dr]|(.7){\Downarrow\lambda'} & 
        {B'}\ar[d]^{f'}\ar[dr]^{v'} \\
        {} & {C'}\ar[ur]^{\ell'}\ar[r]_{g'}\ar[dr]_{w'} & 
        {A'}\ar[dr]|*+<1pt>{\scriptstyle u'} & {B''}\ar[d]^{f''} \\
        {} & {} & {C''}\ar[r]_{g''} & {A''}
      }
    } 
    \mkern5mu = \mkern5mu
    \vcenter{
      \xymatrix@=1.5em{
        {} & {B}\ar[dr]^{v} & {} & {} \\
        {C}\ar[ur]^{\ell}\ar[dr]_{w} & 
        {\scriptstyle\Downarrow\tau} & 
        {B'}\ar[dr]^{v'} \\
        {} & {C'}\ar[ur]^{\ell'}\ar[dr]_{w'} & 
        {\scriptstyle\Downarrow\tau'} 
        \ar@{}[dr]|(.7){\Downarrow\lambda''}
        & {B''}\ar[d]^{f''} \\
        {} & {} & {C''}\ar[ur]^{\ell''}\ar[r]_{g''} & {A''}
      }
    } 
  \end{equation*}
  which shows that the 2-cell induced by the composite of a pair of transformations may be computed as the pasted composite, shown on the right, of the 2-cells induced by each individual transformation. Now if those induced 2-cells are isomorphisms, then so is their pasted composite, from which our desired composition result follows.
\end{obs}

\subsection*{Projective cofibrant weighted limits in \texorpdfstring{$\protect\qCat_{r\infty}^{\protect\pbshape}$}{qCat}}

Our aim is now to prove:

\begin{prop}\label{prop:abslifts-in-limits}
  The quasi-categorically enriched subcategory $\qCat_{r\infty}^\pbshape$ is closed in $\qCat_\infty^\pbshape$ under limits weighted by projective cofibrant weights.
\end{prop}

Before proving Proposition \ref{prop:abslifts-in-limits}, let us explain  how it will provide us with a proof of Theorem~\ref{thm:limits-in-limits}.

\begin{proof}[Proof of Theorem~\ref{thm:limits-in-limits}]
  For any fixed simplicial set $X$, there exists a simplicially enriched functor $F_X\colon\qCat_\infty\to\qCat_\infty^\pbshape$ which carries a quasi-category $A$ to the diagram
  \begin{equation}\label{eq:limit-X-pbshape}
    \xymatrix@=1.5em{
      {} & {A}\ar[d]^c \\
      {A^X}\ar[r]_{\id} & {A^X}
    }
  \end{equation}
  and carries a functor $f\colon A\to B$ to the transformation:
  \begin{equation}\label{eq:limit-X-pbshape-trans}
    \xymatrix@=1.5em{
      {} & {A}\ar[d]_{c}\ar[dr]^{f} & {} \\
      {A^X}\ar[r]^{\id}\ar[dr]_{f^X} & {A^X}
      \ar[dr]|(0.45)*+<1pt>{\scriptstyle f^X} & 
      {B}\ar[d]^{c} \\
      {} & {B^X}\ar[r]_{\id} & {B^X}
    }
  \end{equation}
  By Definition \refI{defn:families.of.diagrams}, the quasi-category $A$ possesses limits of shape $X$ if and only if \eqref{eq:limit-X-pbshape} admits an absolute right lifting, and a functor $f\colon A\to B$ preserves limits of shape $X$ if and only if \eqref{eq:limit-X-pbshape-trans} is right exact. In other words, the simplicial category $\qCat_{X\infty}$ of Definition \ref{defn:qCat-X}  is the inverse image of the simplicial subcategory $\qCat_{r\infty}^\pbshape$ along the simplicial functor $F_X\colon\qCat_\infty\to\qCat_\infty^\pbshape$.

  Now $\qCat_\infty^\pbshape$ inherits limits pointwise from $\qCat_\infty$ and the exponentiation functor $(-)^X$ is right adjoint, so  $F_X$ preserves those limits under which $\qCat_\infty\subset\sSet$ is closed. In particular, it preserves all limits weighted by projective cofibrant weights.  It is a general and easily demonstrated fact that if $\tcat{A}$ and $\tcat{B}$ are enriched categories admitting a certain class of limits, $F\colon\tcat{A}\to\tcat{B}$ is an enriched functor preserving those limits, and $\tcat{B}'$ is an enriched subcategory of $\tcat{B}$ which is closed under those limits, then the enriched inverse image subcategory $\tcat{A}'\defeq F^{-1}(\tcat{B}')$ is also closed in $\tcat{A}$ under those limits. 

  In particular, Proposition~\ref{prop:abslifts-in-limits} tells us that the simplicial subcategory $\qCat_{r\infty}^\pbshape$ is closed in $\qCat_\infty^\pbshape$ under limits weighted by projective cofibrant weights, and $F_X\colon\qCat_\infty\to\qCat_\infty^\pbshape$ preserves these limits. So it follows that the inverse image $\qCat_{X\infty}$ of $\qCat_{r\infty}^\pbshape$ along $F_X$ is closed in $\qCat_\infty$ under projective cofibrant weighted limits.
\end{proof}

Now to prove Proposition~\ref{prop:abslifts-in-limits}, we first characterise the objects and 0-arrows of $\qCat_{r\infty}^\pbshape$ in terms of the possession and preservation of terminal objects in certain comma quasi-categories.

\begin{obs}\label{obs:comma-simp-functor}
As described in Example~\refII{ex:commaweightedlimit}, the comma quasi-category construction is one example of a limit weighted by a projective cofibrant weight $W\colon\pbshape\to\sSet$ so, in particular, it gives rise to a simplicial functor $\comma=\wlim{W}{-}\colon\qCat_\infty^\pbshape\to\qCat_\infty$. Weighted limits commute with weighted limits so it follows that this comma construction preserves all limits that are constructed pointwise in $\qCat_\infty^\pbshape$.
\end{obs}

The following pair of results combine to supply a pointwise characterisation of right exactness.

{
\renewcommand{\thethm}{\refI{thm:pointwise}}
  \begin{thm}
A functor $g\colon C\to A$ admits an absolute right lifting through the functor $f\colon B\to A$ if and only if for all objects $c$ of $C$ the quasi-category $f\comma gc$ has a terminal object.
  \end{thm}
  \addtocounter{thm}{-1}
}
{
\renewcommand{\thethm}{\refI{cor:pointwise}}
\begin{cor}
  A triangle \eqref{eq:lifting.triangle}  is an absolute right lifting diagram if and only if for each $c \colon \Del^0 \to C$ the restricted triangle
  \[    \xymatrix{ \ar@{}[dr]|(.7){\Downarrow\lambda c} & B \ar[d]^f \\ \Del^0 \ar[r]_{gc} \ar[ur]^{\ell c} & A}\] displays $\ell c$ as an absolute right lifting of $gc$ through $f$.
\end{cor}
  \addtocounter{thm}{-1}
}

\begin{obs}[a pointwise characterisation of right exactness]\label{obs:pointwise-exactness}
Suppose that~\eqref{eq:obj-qCat-pbshape-trans} is a transformation of diagrams which admit absolute right liftings. Then we may pick any absolute right lifting~\eqref{eq:lifting.triangle} of its domain and apply Corollary~\refI{cor:pointwise} to show that for all objects $c$ of $C$ the 2-cell $\lambda c\colon f\ell c\Rightarrow g c$ induces an object which is terminal in $f\comma gc$. Applying that same result to the composite triangle on the left of~\eqref{eq:induced-mate} we find that~\eqref{eq:obj-qCat-pbshape-trans} is a right exact transformation if and only if for all objects $c$ in $C$ the 2-cell $u\lambda c\colon f'v\ell c = uf\ell c\Rightarrow ugc = g'wc$ induces an object which is terminal in $f'\comma g'wc$. 
  
  Now observe that for each object $c$ of $C$ the transformation~\eqref{eq:obj-qCat-pbshape-trans} restricts to define a transformation
  \begin{equation*}
    \xymatrix@=1.5em{
      {} & {B}\ar[d]_{f}\ar[dr]^{v} & {} \\
      {\Del^0}\ar[r]^{gc}\ar@{=}[dr]_{\id} & {A}\ar[dr]|*+<1pt>{\scriptstyle u} & 
      {B'}\ar[d]^{f'} \\
      {} & {\Del^0}\ar[r]_{g'wc} & {A'}
    }
  \end{equation*}
which is mapped to a functor $v\comma_u \id\colon f\comma gc\to f'\comma g'wc$ under the comma construction. This functor maps any object of $f\comma gc$ which is induced by $\lambda c\colon f\ell c\Rightarrow gc$ to an object of $f'\comma g'wc$ which is induced by $u\lambda c\colon f'v\ell c\Rightarrow g'wc$. 
\end{obs}

  Combining and summarising these facts, we get the following lemma.

\begin{lem}\label{lem:right-exact-pointwise}
  The diagram $C \xrightarrow{g} A \xleftarrow{f} B$ is in the subcategory $\qCat_{r\infty}^\pbshape$ if and only if for all objects $c$ of $C$ the comma quasi-category $f\comma gc$ has a terminal object. A transformation~\eqref{eq:obj-qCat-pbshape-trans} between two such diagrams is in the subcategory $\qCat_{r\infty}^\pbshape$ if and only if for all objects $c$ in $C$ the induced functor $v\comma_u\id\colon f\comma gc \to f'\comma g'wc$ preserves terminal objects.
\end{lem}

\begin{proof}
  The first statement is the content of Theorem~\refI{thm:pointwise}. For the second, fix an absolute right lifting as in~\eqref{eq:lifting.triangle} and an object $c$ of $C$. By Corollary \refI{cor:pointwise},  an object induced by $\lambda c\colon f\ell c\Rightarrow gc$ is terminal in $f\comma gc$. The functor $v\comma_u\id\colon f\comma gc \to f'\comma g'wc$ carries this terminal object to an object induced by $u\lambda c\colon f'v\ell c\Rightarrow g'wc$ in $f'\comma g'wc$. Hence, it preserves terminal objects if and only if any object induced by the 2-cell $u\lambda c$ is terminal in $f'\comma g'wc$. As discussed in Observation \ref{obs:pointwise-exactness}, this latter condition holds for all objects $c$ in $C$ if and only if the transformation in~\eqref{eq:obj-qCat-pbshape-trans} is right exact, as required.
\end{proof}

  The characterisation of the objects and 0-arrows of $\qCat_{r\infty}^\pbshape$ provided by Lemma \ref{lem:right-exact-pointwise} enables the proof of the following result:

\begin{lem}\label{lem:conical-closure}
  The subcategory $\qCat_r^\pbshape$ is closed in $\qCat^\pbshape$ under arbitrary small products, pullbacks of pointwise isofibrations, countable composites of pointwise isofibrations, and splitting of idempotents.
\end{lem}

\begin{proof}
  Suppose that we are given a diagram $D\colon \scat{A}\to\qCat_r^\pbshape$ of one of the kinds described in the statement. We shall adopt the notation
  \begin{equation*}
    \vcenter{
    \xymatrix@=1.5em{
      {} & {B^i}\ar[d]^{f^i} \\
      {C^i}\ar[r]_{g^i} & {A^i}
    }}
    \mkern40mu\text{and}\mkern40mu
    \vcenter{
    \xymatrix@=1.5em{
      {} & {B^i}\ar[d]_{f^i}\ar[dr]^{v^\phi} & {} \\
      {C^i}\ar[r]^{g^i}\ar[dr]_{w^\phi} & {A^i}\ar[dr]|-*+<1pt>{\scriptstyle u^\phi} & 
      {B^j}\ar[d]^{f^j} \\
      {} & {C^j}\ar[r]_{g^j} & {A^j}
    }}
  \end{equation*}
  for the diagram and transformation obtained by evaluating $D\colon \scat{A}\to\qCat_r^\pbshape$ at an object $i$ and at an arrow $\phi\colon i\to j$ of $\scat{A}$ respectively. The diagram $D$ possesses a pointwise limit in $\qCat^\pbshape$ and we shall adopt the notation
  \begin{equation*}
    \xymatrix@R=1.5em@C=3em{
      {} & {\lim_{i\in\scat{A}} B^i}\ar[d]_{f}\ar[dr]^{\pi^{b,i}} & {} \\
      {\lim_{i\in\scat{A}} C^i}\ar[r]^{g}\ar[dr]_{\pi^{c,i}} & {\lim_{i\in\scat{A}} A^i}\ar[dr]^(0.45){\pi^{a,i}} & 
      {B^i}\ar[d]^{f^i} \\
      {} & {C^i}\ar[r]_{g^i} & {A^i}
    }
  \end{equation*}
  for the component of its limit cone at the object $i$ of $\scat{A}$.

  Given an object $c$ of $\lim_{i\in\scat{A}} C^i$ we may apply the limit projection at the object $i$ in $\scat{A}$ to obtain an object $c^i\defeq \pi^{c,i}c$ of $C^i$. The limit projections commute with the arrows in the diagram $D$, so if $\phi\colon i\to j$ is an arrow in $\scat{A}$ then  we know that $w^\phi c^i = w^\phi\pi^{c,i}c = \pi^{c,j}c = c^j$. Consequently we have arrows as displayed on the left-hand side
  \begin{equation}\label{eq:conical-limit-cone}
    \vcenter{
    \xymatrix@=1.5em{
      {} & {B^i}\ar[d]_{f^i}\ar[dr]^{v^\phi} & {} \\
      {\Del^0}\ar[r]^{g^ic^i}\ar@{=}[dr]_{\id} & 
      {A^i}\ar[dr]|-*+<2pt>{\scriptstyle u^\phi} & 
      {B^j}\ar[d]^{f^j} \\
      {} & {\Del^0}\ar[r]_{g^jc^j} & {A^j}
    }}
    \mkern40mu\text{and}\mkern40mu
    \vcenter{
    \xymatrix@R=1.5em@C=3em{
      {} & {\lim_{i\in\scat{A}} B^i}\ar[d]_{f}\ar[dr]^{\pi^{b,i}} & {} \\
      {\Del^0}\ar[r]^-{gc}\ar@{=}[dr]_{\id} & {
      \lim_{i\in\scat{A}} A^i}\ar[dr]^(0.45){\pi^{a,i}} & 
      {B^i}\ar[d]^{f^i} \\
      {} & {\Del^0}\ar[r]_{g^ic^i} & {A^i}
    }}
  \end{equation}
   that assemble to give a diagram $D_c\colon\scat{A}\to\qCat^{\pbshape}$ and a cone over that diagram, as displayed on the right. As limits in $\qCat^\pbshape$ are constructed pointwise in $\qCat$, this cone is a limit cone. Because the limit of any diagram of terminal objects is a terminal object, the object in the lower left-hand corner of the limit diagram is $\Delta^0$.
 
  Applying the comma construction $\comma\colon\qCat^{\pbshape}\to\qCat$, we obtain a diagram $\scat{A}\to\qCat$ which maps each object $i$ of $\scat{A}$ to $f^i\comma g^ic^i$ and maps an arrow $\phi\colon i\to j$ of $\scat{A}$ to the induced functor $v^\phi\comma_{u^\phi}\id\colon f^i\comma g^ic^i\to f^j\comma g^jc^j$. The comma construction, as a weighted limit, preserves limits so it maps the limit cone \eqref{eq:conical-limit-cone} to a limit cone in $\qCat$, which displays $f\comma gc$ as a limit of the diagram $\stcat{A}\to\qCat$.

  Now $D$ is a diagram which lands in $\qCat_r^\pbshape$, mapping objects to diagrams which admit an absolute right lifting and arrows to right exact transformations. So Lemma~\ref{lem:right-exact-pointwise} tells us that for each object $i$ of $\scat{A}$ the quasi-category $f^i\comma g^ic^i$ has a terminal object and that for each arrow $\phi\colon i\to j$ of $\scat{A}$ the functor $v^\phi\comma_{u^\phi}\id\colon f^i\comma g^ic^i\to f^j\comma g^jc^j$ preserves terminal objects. Furthermore, if an arrow $\phi\colon i\to j$ of $\scat{A}$ is mapped to a pointwise isofibration under $D$ then it is also mapped to a pointwise isofibration under $D_c$. Observation~\refI{obs:comma-obj-maps} demonstrates that the comma quasi-category functor carries pointwise isofibrations to isofibrations, so it follows that the functor $v^\phi\comma_{u^\phi}\id\colon f^i\comma g^ic^i\to f^j\comma g^jc^j$ is an isofibration. Summarising these facts we find that the composite functor $\stcat{A} \xrightarrow{D_c} \qCat^\pbshape \xrightarrow{\comma} \qCat$ lands in the subcategory $\qCat_\emptyset$ and parametrises one of the conical limits listed in Observation~\ref{obs:terminals-in-limits-prf}.

  Now depending on the particular limit involved, one of Lemmas~\ref{lem:qCat_t-closed-products}, \ref{lem:qCat_t-closed-pullbacks}, \ref{lem:qCat_t-closed-countable-composite}, or~\ref{lem:qCat_t-closed-idempotents} applies to show that the limit $f\comma gc$ of the diagram $\comma \circ D_c$ possesses a terminal object and that each of the limit cone projections $\pi^{b,i}\comma_{\pi^{a,i}}\id\colon f\comma gc\to f^i\comma g^ic^i$ preserves terminal objects. Consequently, Lemma~\ref{lem:right-exact-pointwise} implies that the limit of the diagram $D$ and each of its limit cone projections are in $\qCat_r^\pbshape$.

  Finally, if the transformation on the left of the following display
  \begin{equation*}
    \vcenter{
    \xymatrix@R=1.5em@C=3em{
      {} & {B'}\ar[d]_{f'}\ar[dr]^{v} & {} \\
      {C'}\ar[r]^{g'}\ar[dr]_{w} & 
      {A'}\ar[dr]|-*+<3pt>{\scriptstyle u} & 
      {\lim_{i\in\scat{A}} B^i}\ar[d]^{f} \\
      {} & {\lim_{i\in\scat{A}} C^i}\ar[r]_{g} & 
      {\lim_{i\in\scat{A}} A^i}
    }}
    \mkern80mu
    \vcenter{
    \xymatrix@=1.5em{
      {} & {B'}\ar[d]_{f'}\ar[dr]^{v^i} & {} \\
      {C'}\ar[r]^{g'}\ar[dr]_{w^i} & 
      {A'}\ar[dr]|-*+<3pt>{\scriptstyle u^i} & 
      {B^i}\ar[d]^{f^i} \\
      {} & {C^i}\ar[r]_{g^i} & 
      {A^i}
    }}
  \end{equation*}
  is induced by the cone of right exact transformations on the right, we must show that it too is right exact. To that end, suppose that $c'$ is an object of $C'$ and recall from the earlier argument that $f\comma gwc'$ is a limit of the diagram $ \comma \circ D_{wc'} \colon \scat{A} \to \qCat_\emptyset$. Furthermore,  the functor $v\comma_u\id\colon f'\comma g'c'\to f\comma gwc'$ is  induced by applying the universal property of that limit to the cone of functors $v^i\comma_{u^i} \id\colon f'\comma g'c'\to f^i\comma g^iw^ic'$. Each functor in the limit cone is induced by a right exact transformation and so  preserves terminal objects by Lemma~\ref{lem:right-exact-pointwise}. Hence, by appeal to the appropriate closure lemma cited above,  this induced functor $v\comma_u\id\colon f'\comma g'c'\to f\comma gwc'$ also preserves terminal objects. Lemma~\ref{lem:right-exact-pointwise} now implies that the induced transformation in the display above is indeed right exact, as required.
\end{proof}

We complete the proof of Proposition \ref{prop:abslifts-in-limits} with the following lemma:

\begin{lem}\label{lem:cotensor-closure}
    The simplicial subcategory $\qCat_{r\infty}^\pbshape$ is closed in $\qCat_\infty^\pbshape$ under cotensoring by an arbitrary simplicial set $X$. Moreover, the transformation induced by cotensoring with a monomorphism $X\inc Y$ is a right exact pointwise isofibration.
\end{lem}

\begin{proof}
  Suppose that the diagram~\eqref{eq:obj-qCat-pbshape} is in $\qCat_{r\infty}^\pbshape$ and has the absolute right lifting~\eqref{eq:lifting.triangle}. Then its cotensor by $X$ in $\qCat_\infty^\pbshape$ is determined pointwise, viz
  \begin{equation}\label{eq:cotensor-in-qCat-pb}
    \xymatrix@=1.5em{ & B^X \ar[d]^{f^X} \\ C^X \ar[r]_{g^X} & A^X}
  \end{equation}
  and it is a matter of a routine argument using properties of the 2-adjunction $-\times X\dashv (-)^X$ on the 2-category $\qCat_2$, as discussed in Observation~\refI{obs:transpose-abs-lifting}, to show that the triangle
  \begin{equation*}
    \xymatrix@=2em{ \ar@{}[dr]|(.7){\Downarrow\lambda^X} & B^X \ar[d]^{f^X} \\ C^X \ar[r]_{g^X} \ar[ur]^{\ell^X} & A^X}
  \end{equation*}
  is also an absolute right lifting. This shows that the cotensor~\eqref{eq:cotensor-in-qCat-pb} is an object of $\qCat_{r\infty}^\pbshape$. 

  The projection transformations of this cotensor are indexed by the vertices $x$ of $X$, with the projection given by evaluation at that vertex (precomposition by $x\colon\Del^0\to X$). The 2-functoriality properties of cotensoring imply that the following pasting equality holds:
  \begin{equation*}
    \vcenter{
      \xymatrix@=1.5em{
        {} \ar@{}[dr]|(.7){\Downarrow\lambda^X}
        & {B^X}\ar[d]^(0.55){f^X}\ar[dr]^{B^x} & {} \\
        {C^X}\ar[ur]^{\ell^X}\ar[r]_(0.6){g^X}\ar[dr]_{C^x} & 
        {A^X}\ar[dr]|*+<1pt>{\scriptstyle A^x} & 
        {B}\ar[d]^{f} \\
        {} & {C}\ar[r]_{g} & {A}
      }
    } 
    \mkern40mu = \mkern40mu
    \vcenter{
      \xymatrix@=1.5em{
        {} & {B^X}\ar[dr]^{B^x} & {} \\
        {C^X}\ar[ur]^{\ell^X}\ar[dr]_{C^x} & 
        \ar@{}[dr]|(.7){\Downarrow\lambda} & 
        {B}\ar[d]^{f} \\
        {} & {C}\ar[ur]^{\ell}\ar[r]_{g} & {A}
      }
    } 
  \end{equation*}
  This demonstrates that the comparison 2-cell induced as in~\eqref{eq:induced-mate} is an identity and thus that the projection indexed by $x$ is right exact.

  All that remains is to show that if we are given a transformation
  \begin{equation}\label{eq:induced-into-cotensor}
    \xymatrix@=1.5em{
      {} & {B'}\ar[d]_{f'}\ar[dr]^{v} & {} \\
      {C'}\ar[r]^{g'}\ar[dr]_{w} & 
      {A'}\ar[dr]|*+<1pt>{\scriptstyle u} & 
      {B^X}\ar[d]^{f^X} \\
      {} & {C^X}\ar[r]_{g^X} & {A^X}
    }
  \end{equation}
  whose composite with each of the projection transformations discussed above is right exact, then that transformation itself is right exact. Write $\tau\colon v\ell\Rightarrow \ell^X w$ for the comparison 2-cell induced as in~\eqref{eq:induced-mate} from \eqref{eq:induced-into-cotensor}. Observation~\ref{obs:comparison-2cell-composite} tells us that the comparison 2-cell associated with the composite of \eqref{eq:induced-into-cotensor} with the projection indexed by a vertex $x$ of $X$ is given by the pasting:
  \begin{equation*}
    \vcenter{
      \xymatrix@=2em{
        {B}\ar[r]^v\ar@{}[dr]|{\Downarrow\tau} & 
        {B^X}\ar[r]^{B^x} & {B} \\
        {C}\ar[u]^\ell\ar[r]_w & {C^X}\ar[u]_{\ell^X}
        \ar[r]_{C^x} & {C}\ar[u]_\ell
      }
    } 
  \end{equation*}
  Consequently, the hypothesis that the transformation \eqref{eq:induced-into-cotensor} composes with each projection to give a right exact transformation asserts that the whiskered 2-cell 
  \begin{equation*}
    \xymatrix@C=3em{{C}\ar@/^1.5ex/[r]!L(0.5)^{v\ell}
      \ar@/_1.5ex/[r]!L(0.5)_{\ell^Xw}\ar@{}[r]|{\Downarrow\tau} &
      {B^X} \ar[r]^{B^x} & {B}}
  \end{equation*}
  is an isomorphism for all vertices $x$ in $X$. Taking transposes under the 2-adjunctions $-\times X\dashv (-)^X$ and $-\times C\dashv (-)^C$ on the 2-category of all simplicial sets (discussed in~\refI{obs:2-cat.of.all.simpsets}), this 2-cell corresponds to the 2-cell
  \begin{equation*}
    \xymatrix@C=3em{{\Del^0}\ar[r]^{x} &
      {X}\ar@/^1.5ex/[r]!L(0.5)\ar@/_1.5ex/[r]!L(0.5)
      \ar@{}[r]|{\Downarrow\hat\tau} &
      {B^C}}
  \end{equation*}
  Now these transpositions preserve the isomorphism property of 2-cells, so this latter 2-cell is an isomorphism for all vertices $x$ of $X$.  It follows, by Observation~\refI{obs:pointwise-iso-reprise}, that $\hat\tau$ is an isomorphism. Consequently, its transpose $\tau$ is also an isomorphism, which completes our proof that the transformation~\eqref{eq:induced-into-cotensor} is right exact as required.
  
Finally, the 2-adjunction $-\times X \dashv (-)^X$ is in fact a parametrised 2-adjunction. A map $i \colon X\inc Y$ induces a 2-natural transformation $(-)^i\colon (-)^Y \to (-)^X$ whose components are isofibrations.   By 2-naturality of $(-)^i$, the following pasting equality holds:
  \begin{equation*}
    \vcenter{
      \xymatrix@=1.5em{
        {} \ar@{}[dr]|(.7){\Downarrow\lambda^Y}
        & {B^Y}\ar[d]^(0.55){f^Y}\ar[dr]^{B^i} & {} \\
        {C^Y}\ar[ur]^{\ell^Y}\ar[r]_(0.6){g^Y}\ar[dr]_{C^i} & 
        {A^Y}\ar[dr]|*+<1pt>{\scriptstyle A^i} & 
        {B^X}\ar[d]^{f^X} \\
        {} & {C^X}\ar[r]_{g^X} & {A^X}
      }
    } 
    \mkern40mu = \mkern40mu
    \vcenter{
      \xymatrix@=1.5em{
        {} & {B^Y}\ar[dr]^{B^i} & {} \\
        {C^Y}\ar[ur]^{\ell^Y}\ar[dr]_{C^i} & 
        \ar@{}[dr]|(.7){\Downarrow\lambda^X} & 
        {B^X}\ar[d]^{f^X} \\
        {} & {C^X}\ar[ur]^{\ell^X}\ar[r]_{g^X} & {A^X}
      }
    } 
  \end{equation*}
  which tells us that the transformation induced by $(-)^i$ is right exact.
\end{proof}

\begin{cor}\label{cor:terminal-cotensor} The category $\qCat_{\emptyset\infty}$ is closed in $\qCat_{\infty}$ under cotensoring by an arbitrary simplicial set $X$. Moreover, the cotensor with an inclusion $X\inc Y$ induces a terminal-object preserving isofibration of quasi-categories.
\end{cor}
\begin{proof}
Specialising Lemma \ref{lem:cotensor-closure} to a diagram
  \begin{equation*}
    \xymatrix@=2em{ & A \ar[d]^{!} \\ \catone \ar@{=}[r]  & \catone}
  \end{equation*}
  in $\qCat_r$ constructs cotensors for a quasi-category $A$ with a terminal object.  The top right-hand component of the right exact pointwise isofibration induced by cotensoring with an inclusion $X\inc Y$ is a terminal-object-preserving isofibration of quasi-categories.
\end{proof}

We now complete the proof of Proposition \ref{prop:abslifts-in-limits} and, in doing so, finish proving Theorem~\ref{thm:limits-in-limits}.

\begin{proof}[Proof of Proposition~\ref{prop:abslifts-in-limits}]
Lemmas \ref{lem:conical-closure} and \ref{lem:cotensor-closure} demonstrate that $\qCat_{r\infty}^\pbshape$ is closed in $\qCat_\infty^\pbshape$ under all of the limit types used in the proof of Proposition~\refII{prop:projwlims}. Consequently, it follows, by the argument given in \S\ref{sec:weighted}, that $\qCat_{r\infty}^\pbshape$ is closed in $\qCat_\infty^\pbshape$ under all limits weighted by projective cofibrant weights.
\end{proof}

\subsection*{Homotopy limits of stable quasi-categories and exact functors}

We close this section with a quick application, proving a theorem advertised in the introduction. Compare with \cite[1.1.4.4]{Lurie:2012uq}.

\begin{thm}\label{thm:stable} 
The simplicial subcategory of stable quasi-categories and exact functors is closed in $\qCat_\infty$ under projective cofibrant weighted limits. In particular, the homotopy limit of a diagram of stable quasi-categories and exact functors is a stable quasi-category.
\end{thm}
\begin{proof}
Theorem \ref{thm:limits-in-limits} implies that a projective cofibrant weighted limit of a diagram of stable quasi-categories and exact functors is finitely complete and cocomplete. It remains to prove that the terminal object is initial and that pullbacks and pushouts coincide.

As demonstrated in the proof of Proposition \refII{prop:projwlims}, a projective cofibrant weighted limit may be constructed as the limit of a conical diagram (an idempotent splitting of the limit of a tower of pullbacks of products of isofibrations) that is built from the original diagram and from cotensors of its objects with monomorphisms of simplicial sets. By the stability hypothesis, in each of the quasi-categories in the original diagram, the terminal object is also initial. The proof of Corollary \ref{cor:terminal-cotensor} constructs the terminal object in a cotensor as the constant functor at the given terminal object; hence terminal objects in cotensors of stable quasi-categories  are also initial. In the proof of Theorem \ref{thm:terminal-closure}, the terminal object in the projective cofibrant weighted limit is constructed as a pointwise terminal object defined with respect to the conical diagram. These pointwise terminals are pointwise initial and hence also initial by the proofs of Lemmas \ref{lem:qCat_t-closed-products}, \ref{lem:qCat_t-closed-pullbacks}, \ref{lem:qCat_t-closed-countable-composite}, and \ref{lem:qCat_t-closed-idempotents}.

The proof that pullbacks coincide with pushouts is similar. We first argue that if pushouts coincide with pullbacks in a quasi-category $A$, then the same is true in $A^X$. A square $\lambda \colon \Delta^1 \times \Del^1 \to A^X$ induces a pair of 2-cells:
\[ \xymatrix{ \ar@{}[dr]|(.7){\Downarrow \lambda} & A^X \ar[d]^c \\ \Delta^0 \ar[ur]^{\lambda_{0,0}} \ar[r]_-{\lambda_{\pbshape}} & (A^X)^\pbshape} \qquad \qquad \xymatrix{ \ar@{}[dr]|(.7){\Uparrow \lambda} & A^X \ar[d]^c \\ \Delta^0 \ar[ur]^{\lambda_{1,1}} \ar[r]_-{\lambda_{\poshape}} & (A^X)^\poshape}\] We use subscripts to denote the evident faces of $\lambda$. As argued in Lemma \ref{lem:cotensor-closure} (see also Proposition \refI{prop:pointwise-limits-in-functor-quasi-categories}), these  define absolute lifting diagrams if and only if the transposed 2-cells do so:
\[ \xymatrix{ \ar@{}[dr]|(.7){\Downarrow \lambda} & A \ar[d]^c \\ X \times \Delta^0 \ar[ur]^{\lambda_{0,0}} \ar[r]_-{\lambda_{\pbshape}} & A^\pbshape} \qquad \qquad \xymatrix{ \ar@{}[dr]|(.7){\Uparrow \lambda} & A \ar[d]^c \\ X \times \Delta^0 \ar[ur]^{\lambda_{1,1}} \ar[r]_-{\lambda_{\poshape}} & A^\poshape}\] 
By Theorem \refI{thm:pointwise} and Corollary \refI{cor:pointwise}, this is the case if and only if this is true upon restricting along each vertex $x \colon \Delta^0 \to X$ (limits and colimits in cotensors are defined pointwise). The data of $\lambda_x$ is a pushout in $A$ if and only if it is a pullback, so we conclude that pullbacks and pushouts in $A^X$ coincide. 

It remains only to consider conical projective cofibrant weighted limits. Let $D \colon \stcat{A} \to \qCat$ be a conical diagram of one of the four types considered by  Lemma \ref{lem:conical-closure} and fix a square $\lambda \colon \Del^1 \times \Del^1 \to \lim_{a \in \stcat{A}} Da$. The proof of Lemma \ref{lem:conical-closure} demonstrates that a 2-cell on the left-hand side
\[ \vcenter{\xymatrix{ \ar@{}[dr]|(.7){\Downarrow \lambda} & \lim_{a \in \stcat{A}} Da \ar[d]^c \\ \Delta^0 \ar[r]_-{\lambda_{\pbshape}} \ar[ur]^{\lambda_{0,0}} & (\lim_{a \in \stcat{A}} Da)^\pbshape \cong \lim_{a \in \stcat{A}} Da^\pbshape}}\qquad
\vcenter{
    \xymatrix@R=1.5em@C=3em{
      {\ar@{}[dr]|(.7){\Downarrow \lambda} } & {\lim_{a\in\scat{A}} Da}\ar[d]^{c}\ar[dr]^{\pi^{a}} & {} \\
      {\Delta^0}\ar[r]_{\lambda_{\pbshape}}\ar[ur]^{\lambda_{0,0}} \ar@{=}[dr] & {\lim_{a\in\scat{A}} Da^\pbshape}\ar[dr]^(0.45){(\pi^{a})^\pbshape} & 
      {Da}\ar[d]^{c} \\
      {} & {\Del^0}\ar[r]_{\pi^a \lambda_{\pbshape}} & {Da^\pbshape}
    }}
\] 
defines an absolute right lifting diagram if and only if each of the projected 2-cells on the right define absolute right lifting diagrams. If the quasi-categories $A_i$ are stable this is the case if and only if this same data defines absolute left lifting diagrams as displayed below-right.
\[ \vcenter{\xymatrix{ \ar@{}[dr]|(.7){\Uparrow \lambda} & \lim_{a \in \stcat{A}} Da \ar[d]^c \\ \Delta^0 \ar[r]_-{\lambda_\poshape} \ar[ur]^{\lambda_{1,1}} & (\lim_{a \in \stcat{A}} Da)^\poshape \cong \lim_{a \in \stcat{A}} Da^\poshape}}\qquad
\vcenter{
    \xymatrix@R=1.5em@C=3em{
      {\ar@{}[dr]|(.7){\Uparrow \lambda} } & {\lim_{a\in\scat{A}} Da}\ar[d]^{c}\ar[dr]^{\pi^{a}} & {} \\
      {\Delta^0}\ar[r]_{\lambda_{\poshape}}\ar[ur]^{\lambda_{1,1}} \ar@{=}[dr] & {\lim_{a\in\scat{A}} Da^\poshape}\ar[dr]^(0.45){(\pi^{a})^\poshape} & 
      {Da}\ar[d]^{c} \\
      {} & {\Del^0}\ar[r]_{\pi^a \lambda_{\poshape}} & {Da^\pbshape}
    }}
\] 
This implies, again by the proof of Lemma \ref{lem:conical-closure}, that the diagram on the left-hand side is an absolute left lifting. Hence, pushouts and pullbacks in $\lim_{a \in \stcat{A}} Da$ coincide.
\end{proof}

\section{Limits and colimits in the quasi-category of algebras}\label{sec:algebras}

The quasi-category of algebras for a homotopy coherent monad on a quasi-category $A$ is defined via a projective cofibrant weighted limit. After recalling the precise details of this construction, we observe an immediate corollary of Theorem \ref{thm:limits-in-limits}: any limit or colimit possessed by $A$ and preserved by the monad is created in the quasi-category of algebras. A special case of this result, appearing as Corollary \ref{cor:monadic-creates-u-split}, is required by %HERE: "relevant to" became "required by"
the quasi-categorical monadicity theorem. The aim of the rest of this section is to prove Theorem \ref{thm:monadic-completeness}, which drops the preservation hypothesis in the case of limits in $A$, in direct analogy with the classical categorical result.

\begin{rec}[homotopy coherent monads and their algebras]\label{rec:algebras}
  Recall, from~\S\refII{sec:formal}, that $\Mnd$ denotes the full simplicial subcategory of the free homotopy coherent adjunction $\Adj$ determined by its object $+$. Its single hom-space $\Mnd(+,+)$ is the category $\Del+$ of finite ordinals and order-preserving maps, which we think of as a simplicial set by identifying it with its nerve (as always). The endo-composition on this hom-space is the ordinal sum functor $\oplus\colon \Del+\times\Del+\to\Del+$ and its identity is the object $[-1]$. Consequently, any simplicial functor $X\colon \Mnd\to\sSet$ may be specified by giving a simplicial set $X$, the image of the object $+\in\Mnd$, equipped with a left action $\cdot\colon\Del+\times X\to X$ of the simplicial monoid $(\Del+,\oplus,[-1])$. Natural transformations of such are simply equivariant maps, that is, simplicial maps which preserve the left action of $\Del+$.

  A {\em homotopy coherent monad\/} on a quasi-category $A$ is a simplicial functor $T\colon\Mnd\to\qCat_\infty$ which maps the object $+$ to $A$. Equivalently, such structures correspond to left actions $\cdot\colon \Del+\times A\to A$ of the simplicial monoid $(\Del+,\oplus,[-1])$ on $A$. The endomorphism $t \defeq [0] \cdot - \colon A \to A$ descends to a monad on the homotopy category of $A$.
  
  The \emph{quasi-category of algebras} $A[t]$ for a homotopy coherent monad is defined to be the limit of $T\colon\Mnd\to\qCat_\infty$ weighted by the simplicial functor $W_- \colon \Mnd \to \sSet$ which maps $+$ to $\Del[t]$ and whose left action by $\Del+$ is again the ordinal sum operation $\oplus\colon\Del+\times\Del[t]\to\Del[t]$; here $\Del[t] \subset \Del$ is the subcategory containing only those morphisms that preserve the top element in each ordinal. That this construction delivers us a quasi-category of algebras, rather than just a mere simplicial set of such, is a consequence of Proposition \refII{prop:projwlims} and the fact that the weight $W_-$ is a projective cell complex (see Lemma~\refII{lem:projcof1}). This quasi-category of algebras comes equipped with a forgetful functor $u^t \colon A[t] \to A$, which is the component of its defining limit cone at $[0]\in\Del[t]$ (see Definition~\refII{ex:monadicadj}).
\end{rec}

\begin{obs}\label{obs:algcats-functoriality}
  The weighted limit construction $\wlim{W_-}{-}\colon\qCat_\infty^{\Mnd}\to\qCat_\infty$ is functorial, so it maps a natural transformation $f\colon A\to B$ of homotopy coherent monads to a functor $f[t]\colon A[t]\to B[t]$ of quasi-categories of algebras. The family of underlying functors are natural with respect to this simplicial functor structure, so in particular the following square
  \begin{equation}\label{eq:commutes-with-underlying}
    \xymatrix@=1.5em{
      {A[t]}\ar@{->>}[d]_{u^t}\ar[r]^{f[t]} &
      {B[t]}\ar@{->>}[d]^{u^t} \\
      {A}\ar[r]_f & {B}
    }
  \end{equation}
  commutes. The simplicial category of homotopy coherent monads $\qCat_\infty^{\Mnd}$ possesses all limits weighted by projective cofibrant weights, because it is a functor category. 
  These are computed as in $\qCat_\infty$, and the quasi-category of algebras simplicial functor $-[t]\defeq\wlim{W_-}{-}\colon\qCat_\infty^{\Mnd}\to\qCat_\infty$ preserves them, because limits commute.
\end{obs}

\begin{obs}
  Every $0$-arrow of $\Mnd$ is an iterated composite of the 0-arrow $[0]\in\Mnd(+,+)$, so it follows that all of the functors in the image of a homotopy coherent monad $T\colon\Mnd\to\qCat_\infty$ preserve some limit or colimit in the quasi-category $A$ if and only if the single functor $t\colon A\to A$ to which $[0]\in\Mnd(+,+)$ maps under $T$ preserves that limit or colimit. We call $t\colon A\to A$ the {\em functor part\/} of the homotopy coherent monad $T$.
\end{obs}

Now Theorem~\ref{thm:limits-in-limits} admits the following immediate corollary:

\begin{cor}\label{cor:monadic-colimits}
  Let $T\colon\Mnd\to\qCat_\infty$ be a homotopy coherent monad on a quasi-category $A$ and let $X$ be a simplicial set. Suppose also that $A$ admits and the functor part $t\colon A\to A$ of $T$ preserves all (co)limits of shape $X$. Then the forgetful functor $u^t\colon A[t]\to A$ creates all (co)limits of shape $X$.\qed
\end{cor}

Using Corollary \ref{cor:monadic-colimits}, %HERE: changed "we can now give a direct proof" to "it follows immediately" and ssec:monadicity to sec:monadicity
it follows immediately that the monadic adjunction satisfies one of the key hypotheses in the quasi-categorical monadicity theorem of \S\refII{sec:monadicity}.

\begin{cor}\label{cor:monadic-creates-u-split} The monadic forgetful functor $u^t \colon A[t] \to A$ creates colimits of $u^t$-split simplicial objects.
\end{cor}
\begin{proof}
A $u^t$-split simplicial object in $A[t]$ is a diagram of the form
\[ \xymatrix{ \Del\op \ar@{u(->}[d] \ar[r] & A[t] \ar[d]^{u^t} \\ \Del[t] \ar[r] & A}\] In particular, a $u^t$-split simplicial object in $A[t]$ defines a split simplicial object in $A$. Theorem \refI{thm:splitgeorealizations} proves that split simplicial objects are absolute colimits, preserved by any functor, and, in particular, by $t$. Hence, Corollary \ref{cor:monadic-colimits} applies.
\end{proof}

In the case of colimits, Corollary \ref{cor:monadic-colimits} is a direct analogue of the usual categorical result. However, in the case of limits we expect the monadic forgetful functor to create all limits that $A$ admits, regardless of whether or not they are preserved by the functor part of the monad. Our aim in this section is to prove the analogous quasi-categorical result, viz:

\begin{thm}\label{thm:monadic-completeness}
Let $T\colon\Mnd\to\qCat_\infty$ define a homotopy coherent monad on a quasi-category $A$. Then $u^t \colon A[t] \to A$ creates any limits that $A$ admits.
\end{thm}

Our proof of this result, which takes up the remainder of this work, parallels that of Theorem \ref{thm:limits-in-limits}: we first dispense with the special case of terminal objects and then extend this result to the general case of limits of any shape. The main technical tool that we use in both cases is an explicit description of the monadic forgetful functor $u^t \colon A[t] \to A$ as the limit of a tower of isofibrations. This description is formal rather than ad-hoc: it arises from the fact that the weights for these limits define a relative projective cell complex.

\subsection*{Understanding the monadic forgetful functor}

For the remainder of this paper we fix a homotopy coherent monad $T \colon \Mnd \to \qCat_\infty$ on a quasi-category $A$.

\begin{obs}[defining the monadic forgetful functor]
 Let $W_+\colon\Mnd\to\sSet$ denote the unique representable weight, which maps $+$ to $\Del+$ and whose left action by $\Del+$ is the ordinal sum $\oplus\colon\Del+\times\Del+\to\Del+$. Applying Yoneda's lemma in the form given in Example~\refII{ex:rep-weights}, we find that the weighted limit $\wlim{W_+}{T}$ is isomorphic to $A$. 
  The equivariant map $-\oplus[0]\colon\Del+\inc\Del[t]$ defines an inclusion $W_+\inc W_-$; this  is the natural transformation which corresponds, under Yoneda's lemma, to the vertex $[0]$ of $W_-(+) = \Del[t]$. 
  
The inclusion $W_+\inc W_-$ induces a map of weighted limits $A[t] \cong \wlim{W_-}{T}\to\wlim{W_+}{T} \cong A$ which is $u^t\colon A[t]\to A$, by definition. 
\end{obs}

\begin{rec}[graphical calculus for $\Adj$]
Our aim is to show that the inclusion $W_+\inc W_-$ is an explicitly presented relative projective cell complex. To characterise the cells, we recall the graphical calculus introduced in \S\refII{sec:generic-adj} to describe the hom-spaces in the simplicial category $\Adj$, including in particular the simplicial sets $\Del+$ and $\Del[t]$ on which the homotopy coherent monads $W_+$ and $W_-$ act. Specialising Definition \refII{def:strict-und-squiggles}, an $n$-simplex in $\Del[t]$ is a \emph{strictly undulating squiggle on $n+1$ lines} that begins at $-$ and ends at $+$ (composition order is from right to left):
\begin{equation*}
\vcenter{\hbox{
\begin{asy}
  int[] vArray = {6,2,5,3,4,0,6,1,3,0};
  Simplex vSimp = Simplex(5,vArray);
  vSimp.draw();
\end{asy}
}} 
\end{equation*}
The $i\th$ face of the $n$-simplex is represented by the squiggle obtained by removing the $i\th$ line and straightening if necessary to preserve the strict undulation. The $n$-simplex is non-degenerate if and only if each of the characters $\{1,2,\ldots, n\}$ appears as an ``interior turn-around point'' listed in the sequence of characters displayed at the bottom; the 5-simplex displayed above is non-degenerate. In this context, we say an $n$-simplex of $\Del[t]$ is \emph{atomic} if it does not contain the symbol + in its interior; the 5-simplex displayed above is not atomic. The non-atomic $n$-simplices are in the image of the atomic $n$-simplices under the action of the left $\Del+$-action.

The image of $\Del+$ under the inclusion $-\oplus[0]\colon\Del+\inc\Del[t]$ consists of those strictly undulating squiggles whose right-most entry but one is $+$, e.g., the degenerate 3-simplex:
\begin{equation*}
\vcenter{\hbox{
\begin{asy}
  int[] vArray = {4,2,4,0,2,1,4,0};
  Simplex vSimp = Simplex(3,vArray);
  vSimp.draw();
\end{asy}
}}
\end{equation*} is in the image of $\Del+$.
\end{rec}

\begin{lem}\label{lem:W-.decomp}
The inclusion $W_+\inc W_-$ is a relative projective cell complex: it may be expressed as a transfinite (indeed countable) composite of pushouts of projective cells $\boundary\Del^m\times W_+\inc \Del^m\times W_+$ in $\sSet^{\Mnd}$.
\end{lem}
\begin{proof}
The proofs of Corollary~\refII{cor:uTconservative} and Proposition~\refII{prop:projcofchar2} provide an explicit presentation of $W_-$ as a colimit of a countable sequence of pushouts of projective cells $\boundary\Del^m \times W_+ \inc \Del^m \times W_+$. Each cell corresponds to an atomic and non-degenerate $m$-simplex in $\Del[t]$ which is not in the image of the inclusion $-\oplus[0]\colon\Del+\inc\Del[t]$. This inclusion is surjective on $0$-simplices, so we will always have $m \geq 1$.

  These atomic simplices may be indexed as a countable sequence $\{\underline{a}_i\}_{i\in\mathbb{N}}$ in such a way that each arrow is preceded by all of those of strictly smaller \emph{width} (the number of characters in the representing string minus one; see Definition~\refII{def:strict-und-squiggles}) and by those of the same width but smaller dimension. Then we can define $W_k$ to be the smallest sub-weight of $W_-$ which contains the atomic 0-simplex $\underline{u} \defeq (-,+) \in\Del[t]$ along with all of the simplices $\underline{a}_i$ with $i<k$. With this indexing convention, $W_0$ is the image of $W_+\inc W_-$ and $W_- = \bigcup_{k\in\mathbb{N}} W_k$.  The order of our sequence was selected to ensure that the simplicial boundary of the simplex $\underline{a}_k$ is contained within $W_k$; we write $m_k$ for its dimension. As argued in the proof of Proposition ~\refII{prop:projcofchar2}, the weight $W_{k+1}$ may be constructed from $W_k$ as a pushout 
  \begin{equation}\label{eq:projective-cell} \xymatrix{ \partial\Del^{m_k} \times W_+ \ar[d] \ar@{u(->}[r] & \Del^{m_k} \times W_+ \ar[d] \\ W_{k} \ar@{u(->}[r]^{\subset} & W_{k+1} \poexcursion}\end{equation}
  in which the map $\Del^{m_k} \times W_+\to W_{k+1}$ on the right corresponds to the $m_k$-simplex $\underline{a}_k\in W_{k+1}$ by Yoneda's lemma.
\end{proof}

In particular it follows that the functor $u^t\colon A[t]\tfib A$ is an isofibration of quasi-categories.

\begin{cor}\label{cor:W-.decomp}
The functor $u^t \colon A[t] \tfib A$ is the limit of a tower of isofibrations of quasi-categories defined as projective weighted limits.
\end{cor}
\begin{proof}
The contravariant weighted limit functor $\wlim{-}{T}$ carries colimits in $\sSet^{\Mnd}$ to limits in $\sSet$. So it carries the sequence $W_0\inc W_1\inc...\inc W_k\inc...$ of Lemma~\ref{lem:W-.decomp} to a tower of functors 
  \begin{equation}\label{eq:tower} A[t]\cong \wlim{W_-}{T}\to\cdots \to \wlim{W_{k+1}}{T} \to \wlim{W_{k}}{T} \to \cdots \to \wlim{W_0}{T} \cong A\end{equation}
whose countable composite is $u^t\colon A[t]\to A$. Furthermore, it carries each pushout~\eqref{eq:projective-cell} to a pullback
  \begin{equation}\label{eq:tower.step.pb}
    \xymatrix{
      \wlim{W_{k+1}}{T} \pbexcursion \ar[r] \ar@{->>}[d] & 
      \wlim{\Del^{m_k}\times W_+}{T} 
      \ar@{}[r]|-{\textstyle\cong} \ar@{->>}[d] &
      A^{\Del^{m_k}} \ar@{->>}[d]  \\ 
      \wlim{W_{k}}{T} \ar[r] & 
      \wlim{\partial\Del^{m_k}\times W_+}{T} 
      \ar@{}[r]|-{\textstyle\cong} &
      A^{\partial\Del^{m_k}}
    }
  \end{equation}
which describes the $k\th$ step in the tower decomposition of $u^t$ as a pullback of an isofibration $A^{\Del^{m_k}}\tfib A^{\boundary\Del^{m_k}}$. In particular, this shows that each functor in that tower is an isofibration and thus that its countable composite $u^t$ is also an isofibration.
\end{proof}

\begin{obs}[The key fact]\label{obs:W-.decomp.key}
  Observe that $\underline{u}$ is the final vertex of each atomic simplex $\underline{a}_k$. This follows because $\underline{a}_k$ has $+$ at its left-hand end, $-$ at its right-hand end, and does not contain $+$ in its interior, so when drawn as a squiggle it must have the following form:
\begin{equation*}
\vcenter{\hbox{
\begin{asy}
  int[] vArray = {6,2,5,3,4,0,5,1,3,0};
  Simplex vSimp = Simplex(5,vArray);
  vSimp.draw();
\end{asy}
}}
\end{equation*}
  in which it crosses the lowest line exactly once, with that crossing situated at the far left of the squiggle. From this description of the shape of $\underline{a}_k$ it is now easy to apply Observation~\refII{obs:squiggle-vertices} to show that its final vertex is $\underline{u}$ as stated. This tells us that the inclusion $W_+\inc W_k$ factors through the attaching map of \eqref{eq:projective-cell} via the inclusion $\fbv{{m_k}}\colon\Del^0\inc\partial\Del^{m_k}$ of the final vertex: 
\begin{equation}\label{eq:terminal-fact}
\xymatrix{ & \boundary\Del^{m_k} \times W_+ \ar[d] \\ W_+  \ar@{u(->}[r] \ar[ur]^{\fbv{{m_k}}\times W_+} & W_{k}}
\end{equation}

This fact is {\em absolutely key\/} to our proof that no assumptions are required of $t\colon A\to A$ in order for the quasi-category of algebras $A[t]$ to inherit limits from $A$. To build our intuition for why this is the case, we consider a simpler sub-problem, that of showing that a terminal object $a$ of $A$ admits a canonical $T$-algebra structure. We know that $u^t\colon A[t]\tfib A$ is the functor induced on weighted limits by the inclusion $W_+\inc W_-$, so this problem reduces to constructing the dashed lift in the following diagram:
  \begin{equation*}
    \xymatrix@=1.5em{
      {W_+}\ar@{u(->}[d]\ar[r]^{a} & {A} \\
      {W_-}\ar@{-->}[ur]
    }
  \end{equation*}
  of equivariant maps. Here the horizontal map is the one which corresponds to the terminal object $a$ of $A$ by Yoneda's lemma, and which therefore maps the representing element $\underline{u}$ to that object. We construct this extension inductively, filling a simplex boundary at each step to find a simplex in $A$ to which we might map the next member of our sequence $\{\underline{a}_k\}_{k\in\mathbb{N}}$ of atomic and non-degenerate arrows, closing up at each step under the action of $\Del+$. It is the result of the last paragraph which ensures that we can find these fillers at each step: it tells us that $\underline{u}$ is the terminal vertex of each $\underline{a}_k$ and, thus, that its boundary maps to a sphere in $A$ whose terminal vertex is the terminal object $a$, this being the image of $\underline{u}$ under the original map we are extending. Consequently, we can use the universal property of $a$ to find the filler we seek. 

  Notice that the same argument does not deliver us a side-condition-free result which supplies canonical algebra structures for initial objects. Such a result would require us to know that the initial vertex of each $\underline{a}_k$ must map to an initial object of $A$ under the extension we are building. However, all we know is that these map to some object obtained by applying an iterated composite of the functor part $t\colon A\to A$ of our monad to that initial object. This is not, in general, good enough for our purposes here, unless we happen to know that $t$ preserves initial objects.
\end{obs}

\subsection*{Creation of terminal objects}

Our aim in this section is to prove the following result:

\begin{prop}\label{prop:monadic-terminal}
 The functor $u^t \colon A[t] \tfib A$ creates terminal objects. Explicitly, we show that if $a$ is a terminal object in $A$ then any commutative diagram
\begin{equation}\label{eq:monadic-terminal}
\xymatrix{  & \boundary\Del^n \ar[r]^v \ar[d] & A[t] \ar@{->>}[d]^{u^t} \\ \Del^0 \ar[r]^{\fbv{n}} \ar@/_1pc/[rr]_a & \Del^n \ar[r]^w \ar@{-->}[ur] & A}\end{equation}
 of solid arrows admits the dashed diagonal filler shown.
\end{prop}

\begin{rmk} The $n=0$ case of the lifting property asserted by Proposition~\ref{prop:monadic-terminal} implies that terminal vertices lift along $u^t$.  The $n \geq 1$ cases then allow us to apply Lemma~\ref{lem:rel-lift-terminal-conv} to show that any such lifted vertex is terminal in $A[t]$. 
\end{rmk}

  Before proving this result, observe that it immediately entails the following corollary:

\begin{cor}\label{cor:monadic-terminal} 
  Suppose that $f\colon A\to B$ is a natural transformation of monads, that the quasi-categories $A$ and $B$ possess terminal objects, and that the functor $f$ preserves terminal objects. Then the induced functor $f[t]\colon A[t]\to B[t]$ preserves the terminal objects that are guaranteed by Proposition~\ref{prop:monadic-terminal}.
\end{cor}

\begin{proof}
  This is an easy argument involving the commutative square~\eqref{eq:commutes-with-underlying} and the fact that, by Proposition~\ref{prop:monadic-terminal}, the underlying functors there create (and thus both preserve and reflect) terminal objects.
\end{proof}

The results discussed in Corollary \ref{cor:W-.decomp} and Observation \ref{obs:W-.decomp.key} provide all of the tools that we need to prove the lifting result described in the statement of Proposition~\ref{prop:monadic-terminal}. Our argument directly generalises that given at the end of Observation \ref{obs:W-.decomp.key} by giving an inductive construction of a lift along the tower of isofibrations induced by applying the weighted limit functor $\wlim{-}{T}$ to the relative projective cell decomposition of $W_+\inc W_-$.

\begin{proof}[Proof of Proposition \ref{prop:monadic-terminal}]
Suppose we are given a terminal object $a \in A$ and a lifting problem \eqref{eq:monadic-terminal}. We construct a diagonal filler by inductively lifting along the tower of fibrations \eqref{eq:tower}. By the inductive hypothesis, we may assume that we have already constructed a lift $l_k\colon\Del^n\to \wlim{W_k}{T}$, which features in the following commutative diagram
\begin{equation}\label{eq:ind.step.filler.diag}
\xymatrix{ & & A[t] \ar@{->>}[d] \\ & \partial\Del^n \ar[ur]^v \ar[r] \ar@{u(->}[d] & \wlim{W_{k+1}}{T} \pbexcursion \ar[r] \ar@{->>}[d]  & A^{\Del^{m_k}} \ar@{->>}[d]  \\ \Del^0 \ar[r]^-{\fbv{n}} \ar[rrd]_a & \Del^n \ar@{-->}[ur]^(0.45){l_{k+1}}\ar@{-->}@/_/[urr]_(0.65)m \ar[dr]^(0.55)w \ar[r]_(0.45){l_k} & \wlim{W_k}{T} \ar@{->>}[d] \ar[r] & A^{\partial\Del^{m_k}} \ar[dl]^{\ev_{\fbv{{m_k}}}} \\ & & A}  
\end{equation}
  in which the right-hand square is the pullback displayed in~\eqref{eq:tower.step.pb} and the lower right-hand commutative triangle is derived by applying $\wlim{-}{T}$ to the commutative triangle \eqref{eq:terminal-fact}.

  The inductive step of our argument requires us to demonstrate the existence of the dashed lift $l_{k+1}\colon\Del^n\to{\wlim{W_{k+1}}{T}}$ filling the left-hand square in this diagram. By the pullback property of the right-hand square, it suffices to construct a lift $m\colon\Del^n\to A^{\Del^{m_k}}$ for the composite of these two squares. Transposing that latter lifting problem, we find in turn that we must solve a corresponding lifting problem:
  \[ \xymatrix{ \partial\Del^n \times \Del^{m_k} \cup \Del^n\times\partial\Del^{m_k} \ar@{u(->}[d] \ar[r]^-f & A \\ \Del^n \times \Del^{m_k} \ar@{-->}[ur]}\] 
The commutativity of the triangles at the bottom of~\eqref{eq:ind.step.filler.diag} implies that the horizontal map $f$ at the top of this transposed problem carries the final vertex $(\fbv{n},\fbv{{m_k}})$ to the terminal object $a$ of $A$. Moreover, all of the simplices of the prism $\Del^n\times\Del^{m_k}$ that are not contained in its boundary $\boundary\Del^n \times \Del^{m_k} \cup \Del^n\times\boundary\Del^{m_k}$ have final vertex $(\fbv{n},\fbv{{m_k}})$, so it follows that we can use the universal property of the terminal object $a$ to extend $f$ from one skeleton of $\Del^n\times\Del^{m_k}$ to the next in order to construct the desired filler and thereby complete our proof.
\end{proof}

\subsection*{Creation of limits}

We now move on to the general case. Given a diagram $g \colon X \to A[t]$, we wish to show that if $u^t g \colon X \to A$ has a limit, then this lifts to a limit of $g$ in $A[t]$ that is preserved by $u^t \colon A[t] \to A$. That is, our aim will be to construct an absolute right lifting 
\begin{equation}\label{eq:lift1} \xymatrix{ & A[t] \ar[d]^c \\ \catone \ar[r]_-g & A[t]^X}\end{equation} of $g$ along the constant diagram map. It will be most economical to describe the construction of this absolute right lifting diagram in the following slightly more general context.

\begin{obs} A homotopy coherent monad on $A$ induces a pointwise-defined homotopy coherent monad on $A^X$ in such a way that the constant diagram map $c\colon A \to A^X$ is a natural transformation of homotopy coherent monads. As weighted limits commute, the diagram \eqref{eq:lift1} is isomorphic to the diagram:
\[ \xymatrix{ & A[t] \ar[d]^{c[t]} \\ \catone \ar[r]_-g & A^X[t]}\] 
Consequently, Theorem~\ref{thm:monadic-completeness} may be obtained as an immediate corollary to the following proposition, simply by specialising it to the natural transformation $c \colon A \to A^X$.
\end{obs}

Our aim is to prove the following general result.

\begin{thm}\label{thm:monadic-limit}
 Suppose $f \colon B \to A$ is a natural transformation of homotopy coherent monads, and suppose that 
\begin{equation}\label{eq:monadic-hypothesis} \xymatrix{ & \ar@{}[dr]|(.6){\Downarrow\lambda} & B \ar[d]^f \\ \catone \ar[r]_-g \ar[urr]^\ell  & A[t] \ar[r]_-{u^t} & A}\end{equation} is an absolute right lifting diagram. Then 
\begin{equation}\label{eq:monadic-conclusion}
  \xymatrix{ & B[t] \ar[d]^{f[t]} \\ \catone \ar[r]_-g & A[t]}
\end{equation}
admits an absolute right lifting. Furthermore, that lifting is preserved by underlying functors, in the sense that the transformation
  \begin{equation}\label{eq:pres-by-underlying}
    \xymatrix@=1.5em{
      {} & {B[t]}\ar[d]_{f[t]}\ar@{->>}[dr]^{u^t} & {} \\
      {\catone}\ar[r]^-{g}\ar@{=}[dr]_{\id} & {A[t]}\ar@{->>}[dr]|*+<1pt>{\scriptstyle u^t} & 
      {B}\ar[d]^{f} \\
      {} & {\catone}\ar[r]_{u^tg} & {A}
    }
  \end{equation}
is right exact.
\end{thm}

\begin{rmk}
  Observation~\ref{obs:pointwise-exactness} reminds us that the diagrams~\eqref{eq:monadic-hypothesis} and~\eqref{eq:monadic-conclusion} admit absolute right liftings if and only if the comma quasi-categories $f\comma u^tg$ and $f[t]\comma g$ admit terminal objects. Furthermore, Lemma~\ref{lem:right-exact-pointwise} tells us that the isofibration $f[t]\comma g\tfib f\comma u^tg$ induced by the pointwise isofibration~\eqref{eq:pres-by-underlying} preserves these terminal objects if and only if~\eqref{eq:pres-by-underlying} is right exact. These observations reduce Theorem~\ref{thm:monadic-limit} to the following equivalent result:
\end{rmk}

\begin{prop}\label{prop:monadic-limit-2}
  Suppose that the comma quasi-category $f\comma u^tg$ possesses a terminal object. Then the comma quasi-category $f[t]\comma g$ also possesses a terminal object and it is preserved by the isofibration $f[t]\comma g\tfib f\comma u^tg$ induced by~\eqref{eq:pres-by-underlying}.
\end{prop}

  We now proceed to prove this latter proposition, but first we require a few preparatory observations:

\begin{obs}\label{obs:expl-alg-desc}
  The defining universal property of weighted limits tells us that simplicial maps $h\colon X\to A[t]$ correspond to equivariant maps $\bar{h}\colon X\times W_-\to A$, where $X\times W_-$ denotes the tensor of $W_-$ by the simplicial set $X$. These corresponding maps are related by the following commutative diagram:
    \begin{equation} \label{eq:expl-alg-desc}\xymatrix{ X \ar[r]^h \ar[d]_{\id_X\times [0]} & A[t] \ar[d]^{u^t} \\ {X\times W_-} \ar[r]_-{\overline{h}} & A}\end{equation}
  This result also provides an explicit description of the action of the functor $f[t]\colon A[t]\to B[t]$ constructed from a natural transformation $f\colon A\to B$ of homotopy coherent monads. This follows by composing the square above with the commutative square in~\eqref{eq:commutes-with-underlying}, which reveals that the composite $f[t]h\colon X\to B[t]$ corresponds to the equivariant map $f\bar{h}\colon X\times W_-\to B$. 

  Taking the special case of $X=\Del^n$ we see that an $n$-simplex $a$ of $A[t]$ may be identified with an equivariant map $\bar{a}\colon\Del^n \times W_-\to A$ and that $f[t]$ acts on that simplex to carry it to the simplex of $B[t]$  identified with the composite equivariant map $f\bar{a}\colon \Del^n\times W_-\to B$. 
\end{obs}

\begin{obs}\label{obs:W-.monad}
  The weight $W_-$ may equally be thought of as giving us a monad on the category $\Del[t]$. Here we have dropped the adjective ``homotopy coherent'' because a homotopy coherent monad on a category is indeed no more nor less than a classical monad, a fact which follows immediately from the characterisation of $(\Del+,\oplus,[-1])$ as the free monoidal category containing a monoid. It is also easily checked that the quasi-category of algebras for a monad on a category is indeed the usual  Eilenberg-Moore category of algebras for that monad. In the particular case of $W_-$, we may show that the category of algebras $\Del[t][t]$ can be identified with the subcategory of top and bottom preserving maps. 

Via Observation \ref{obs:expl-alg-desc}, the identity map $\id\colon W_-=\Del[t] \to \Del[t]$ corresponds to an object $\overline{[0]}$ of $\Del[t][t]$ which is mapped to the object $[0]$ of $\Del[t]$ by the underlying functor $u^t\colon \Del[t][t]\tfib \Del[t]$. Applying Proposition~\ref{prop:monadic-terminal}, we know that $u^t$ creates terminal objects, and we also know that $[0]$ is terminal in $\Del[t]$, so it follows that $\overline{[0]}$ is terminal in $\Del[t][t]$. Furthermore, if $a$ is an object of $A[t]$ given as a monad transformation $\bar{a}\colon \Del[t]\to A$ then the induced functor $\bar{a}[t]\colon \Del[t][t]\to A[t]$ of quasi-categories of algebras maps the object $\overline{[0]}$ of $\Del[t][t]$ to the object $a$ of $A[t]$. This is an immediate consequence of the fact that $\overline{[0]}$ is defined to be the object corresponding to the identity on $\Del[t]$ and the explicit description of $\bar{a}[t]$ given at the end of Observation~\ref{obs:expl-alg-desc}.
\end{obs}

\begin{obs}
  Returning now to the diagrams in the statement of Theorem~\ref{thm:monadic-limit}, we can apply Observation~\ref{obs:expl-alg-desc} to the map $g\colon \catone\to A[t]$ of~\eqref{eq:monadic-conclusion} to obtain a corresponding equivariant map $\bar{g}\colon W_-\to A$. Applying the quasi-category of algebras construction to the monads on $A$, $B$, and $\Del[t]$ and to the monad transformations $f$ and $\bar{g}$, we obtain the following commutative diagram:
  \begin{equation}\label{eq:pre-commas}
    \xymatrix@R=1.75em@C=2.5em{
      {\catone} \ar[r]_-{\overline{[0]}}\ar@{=}[d] 
      \ar@/^1pc/[rr]^{g} &
      {\Del[t][t]}\ar[r]_-{\bar{g}[t]}
      \ar@{->>}[d]_(0.4){u^t} &
      {A[t]}\ar@{->>}[d]^(0.4){u^t} & 
      {B[t]}\ar@{->>}[d]^(0.4){u^t}\ar[l]_-{f[t]} \\
      {\catone} \ar[r]^-{[0]}
      \ar@/_1pc/[rr]_{u^tg} &
      {\Del[t]} \ar[r]^-{\bar{g}} &
      {A} & {B}\ar[l]^-{f}
    }
  \end{equation}
Here the equality $\bar{g}[0] = u^t g$ is a special case of the commutative diagram~\eqref{eq:expl-alg-desc}, and the equality $(\bar{g}[t])\overline{[0]}=g$  is discussed at the end of Observation~\ref{obs:W-.monad}.

  A routine pullback computation, starting from Definition~\refI{def:comma-obj}, reveals that if we are given three functors $f\colon B\to A$, $g\colon C\to A$, and $h\colon D\to C$ then we may express the comma quasi-category $f\comma gh$ as the following pullback of the comma quasi-category $f\comma g$:
  \begin{equation*}
    \xymatrix@=1.5em{
      {f\comma gh}\pbexcursion\ar[r]\ar@{->>}[d]_{p_1} &
      {f\comma g}\ar@{->>}[d]^{p_1} \\
      {D}\ar[r]_h & {C}
    }
  \end{equation*}
  Here the horizontal arrow at the top of this square is obtained by applying the comma construction functor $\comma\colon\qCat_\infty^\pbshape\to\qCat_\infty$ to the following transformation:
  \begin{equation*}
    \xymatrix@=1.5em{
      {D}\ar[r]^{gh}\ar[d]_{h} & {A}\ar@{=}[d] &
      {B}\ar[l]_f\ar@{=}[d] \\
      {C}\ar[r]_g & {A} & {B}\ar[l]_f
    }
  \end{equation*}

Applying this result to the rows of~\eqref{eq:pre-commas}, we obtain the following diagram of commutative squares:
  \begin{equation}\label{eq:squares-of-power}
    \vcenter{
      \xymatrix@=1.5em{
        {f[t]\comma g}\ar@{}[dr]|{\textstyle\text{(A)}}
        \ar[r]\ar@{->>}[d] &
        {f[t]\comma \bar{g}[t]}\ar@{->>}[d] \\
        {f\comma u^tg}\pbexcursion
        \ar@{}[dr]|{\textstyle\text{(B)}}
        \ar[r]\ar@{->>}[d]_{p_1} &
        {f\comma \bar{g}}\ar@{->>}[d]^{p_1} \\
        {\catone}\ar[r]_{[0]} & {\Del[t]}
      }
    } =
    \vcenter{
      \xymatrix@=1.5em{
        {f[t]\comma g}\pbexcursion
        \ar@{}[dr]|{\textstyle\text{(C)}}
        \ar[r]\ar@{->>}[d]_{p_1} &
        {f[t]\comma \bar{g}[t]}\ar@{->>}[d]^{p_1} \\
        {\catone}\ar[r]_{[0]} & {\Del[t]}
      }
    }
  \end{equation}
  Here the pullbacks labelled (B) and (C) are those obtained by applying the result of the last paragraph to the bottom and top rows of~\eqref{eq:pre-commas} respectively. The vertical isofibrations of the square labelled (A) are those induced by the obvious transformations whose legs are various of the vertical isofibrations in \eqref{eq:pre-commas}. Most crucially, the vertical isofibration on the left of square (A) is precisely the induced functor that features in the statement of Proposition~\ref{prop:monadic-limit-2}.
\end{obs}

  This provides us with almost all of the machinery we need to establish Proposition~\ref{prop:monadic-limit-2}. The last piece of the jigsaw is the following lemma:

\begin{lem}\label{lem:technical-lem}
  Under the assumption that~\eqref{eq:monadic-hypothesis} is an absolute right lifting diagram, the comma quasi-category $f\comma \bar{g}$ possesses a terminal object which is preserved by the projection functor $p_1\colon f\comma\bar{g}\tfib\Del[t]$.
\end{lem}

  We postpone the proof of this technical lemma in our brimming enthusiasm to get to our ultimate goal:

\begin{proof}[Proof of Proposition \ref{prop:monadic-limit-2}]
  Our proof revolves around an analysis of the squares displayed in~\eqref{eq:squares-of-power}. First consider the square labelled (B). The domain and codomain of its lower horizontal functor $[0]\colon \catone\to\Del[t]$ both possess terminal objects and the inclusion preserves them, since $[0]$ is the terminal object of $\Del[t]$. What is more, Lemma~\ref{lem:technical-lem} tells us that the right-hand vertical $f\comma \bar{g}\tfib \Del[t]$ of (B) is also a functor whose domain and codomain possess and which preserves terminal objects. So, applying Lemma~\ref{lem:qCat_t-closed-pullbacks}, we find that the top left vertex $f\comma u^tg$ of (B) also possesses a terminal object and that its projection functors preserve such.

The cancellation lemma for pullbacks tells us that the square labelled (A) is also a pullback. So our intention is again to apply Lemma~\ref{lem:qCat_t-closed-pullbacks} to (A) to show that its top left-hand vertex $f[t]\comma g$ possesses a terminal object and that this is preserved by the left-hand vertical isofibration $f[t]\comma g\tfib f\comma u^tg$, these being precisely the conclusions of Proposition~\ref{prop:monadic-limit-2}. We have just shown that the lower horizontal map of (A) is a functor whose domain and codomain possess and which preserves terminal objects, so all that remains in order to apply Lemma~\ref{lem:qCat_t-closed-pullbacks} and complete our proof is to demonstrate the same properties for its right-hand vertical $f[t]\comma \bar{g}[t]\tfib f\comma \bar{g}$. 

  Now Observation~\ref{obs:algcats-functoriality} points out that $\qCat_\infty^{\Mnd}$ possesses all limits weighted by projective cofibrant weights, that these are constructed as in $\qCat_\infty$, and that the quasi-category of algebras construction preserves them. So, applying that result to the comma construction for the monad transformations $f$ and $\bar{g}$, we find that the comma quasi-category $f\comma\bar{g}$ inherits a monad structure and that there exists a canonical isomorphism between the category of algebras $(f\comma\bar{g})[t]$ associated with that monad and the comma quasi-category $f[t]\comma\bar{g}[t]$. What is more, that isomorphism makes the following triangle commute:
  \begin{equation*}
    \xymatrix@C=1em@R=2em{
      {f[t]\comma\bar{g}[t]}\ar[rr]^{\textstyle\cong} 
      \ar@{->>}[dr] && {(f\comma\bar{g})[t]}
      \ar@{->>}[dl]^{u^t} \\
      & {f\comma\bar{g}} &
    }
  \end{equation*}
  where the diagonal map on the left is the right-hand vertical of square (A). Consequently, the required result is equivalent to showing that the quasi-category of algebras $(f\comma\bar{g})[t]$ has a terminal object which is preserved by the underlying functor down to $f\comma\bar[g]$. These facts follow immediately from Proposition~\ref{prop:monadic-terminal} and the fact that we know, from Lemma~\ref{lem:technical-lem}, that $f\comma\bar{g}$ possesses a terminal object.
\end{proof}

Finally, we complete our argument by establishing the technical Lemma~\ref{lem:technical-lem}, which is a special case of the following result:

\begin{lem}
  Suppose that $c$ is a terminal object in $C$ and that we have an absolute right lifting diagram
  \begin{equation*}
    \xymatrix@=1.5em{
      & \ar@{}[dr]|(0.6){\Downarrow\lambda} & {B}\ar[d]^f \\
      {\catone}\ar[r]_c\ar[urr]^\ell &
      {C}\ar[r]_g & {A}
    }
  \end{equation*}
  Then the comma quasi-category $f\comma g$ admits a terminal object which is preserved by the projection $p_1\colon f\comma g\tfib C$.
\end{lem}

\begin{proof}
  All of the components that we need are contained in the proof of Theorem~\refI{thm:pointwise}. Firstly, the existence of the absolute right lifting in the statement is equivalent to asking for the comma quasi-category $f\comma gc \simeq \slicer{f}{gc}$ to have a terminal object, this being the object of $f\comma gc$ induced by $\lambda$. Given that result, the first part of the proof of Theorem~\refI{thm:pointwise} demonstrates that any lifting problem
  \begin{equation*}
    \xymatrix@=1.5em{
      {\boundary\Del^m}\ar[r]\ar@{u(->}[d] &
      {f\comma g}\ar@{->>}[d]^{p_1} \\
      {\Del^m}\ar@{-->}[ur]\ar[r] & {C}
    } 
  \end{equation*} 
  which maps the terminal vertex of $\boundary\Del^m$ to the object of $f\comma g$ induced by $\lambda$ has a solution as marked. However, that object is mapped to the terminal object $c$ of $C$ by the projection $p_1\colon f\comma g\tfib C$, so it follows that we can apply Lemma~\ref{lem:rel-lift-terminal-conv} to complete our proof.
\end{proof}

The desired Theorem \ref{thm:monadic-completeness} is an immediate corollary of Proposition \ref{prop:monadic-terminal} and Theorem \ref{thm:monadic-limit}, or equivalently, Proposition \ref{prop:monadic-limit-2}.

% actual footer stuff
%!TEX root=all.tex
% ******************************************************************
% ** Title:            Homotopy Coherent Adjunctions
% **                   Standard Footer
% ** Precis:        
% ** Author:           Emily Riehl and Dominic Verity
% ** Commenced:        2/3/2012
% ******************************************************************

\ifnum\value{nestcount}>1\global\addtocounter{nestcount}{-1}\else

% actual footer stuff

    \bibliographystyle{abbrv}
    \bibliography{../common/index}
    
   % \listoffixmes

    \end{document}

\fi